\documentclass[leqno,a4paper,12pt]{article}
\usepackage{amsmath,amsthm}

\newtheorem{Thm}{\indent Theorem}[section]
\newtheorem{Prop}[Thm]{\indent Proposition}

\newtheorem{Cor}[Thm]{\indent Corollary}

\theoremstyle{definition}
\newtheorem{Def}[Thm]{\indent Definition}
\newtheorem{Rem}[Thm]{\indent Remark}
\newtheorem{Ex}[Thm]{\indent Example}
\newtheorem{Ass}[Thm]{\indent Assumption}
\def\qed{{\hskip0pt\unskip\unskip\nobreak\hfil\penalty50
          \hskip1em\hbox{}\nobreak\hfil
          {\bf q.e.d.}%
          \parfillskip=0pt\finalhyphendemerits=0
          \par}\medskip}

\newenvironment{Proof}
               {{\it Proof.}\quad}
               {\qed}

\newenvironment{Proofof}[1]
               {{\it Proof of #1.}\quad}
               {\qed}

\newcommand{\Prime}{\kern3\fontdimen1\font$'$\kern-7\fontdimen1\font}

\long\def\forget#1{}

\long\def\beginSIDEREMARK#1\endSIDEREMARK
    {{\par\bigskip\advance\leftskip by 2cm
                  \advance\rightskip by -2cm\noindent
      {\bf Our own side remark:} #1
      \par\bigskip\noindent}}
\long\def\beginFORGET#1\endFORGET{#1}
\long\def\beginFORGET#1\endFORGET{}

\def\?{\ ???\ \immediate\write16{}%
\immediate\write16{Warning: There was still a question mark . . . }%
\immediate\write16{}}

\usepackage{amsmath}
\usepackage{exscale}
\usepackage{amssymb}

\newcommand{\BC}{{\mathbb{C}}}

\newcommand{\BG}{{\mathbb{G}}}

\newcommand{\BN}{{\mathbb{N}}}

\newcommand{\BQ}{{\mathbb{Q}}}
\newcommand{\BR}{{\mathbb{R}}}

\newcommand{\BZ}{{\mathbb{Z}}}

\newcommand{\Fp}{{\mathfrak{p}}}

\newcommand{\FS}{{\mathfrak{S}}}

\newcommand{\FZ}{{\mathfrak{Z}}}

\newcommand{\CA}{{\cal A}}

\newcommand{\CC}{{\cal C}}
\newcommand{\CD}{{\cal D}}

\newcommand{\CH}{{\cal H}}

\newcommand{\CM}{{\cal M}}

\newcommand{\CW}{{\cal W}}

\newcommand{\Spec}{\mathop{{\bf Spec}}\nolimits}

\newcommand{\reg}{\mathop{{\rm reg}}\nolimits}

\newcommand{\Gr}{\mathop{\rm Gr}\nolimits}
\newcommand{\Hom}{\mathop{\rm Hom}\nolimits}

\newcommand{\loccit}{[loc.$\;$cit.]}

\def\tei{\, | \,}
\def\halb{\frac{1}{2}}

\newbox\mybox
\def\arrover#1{\mathrel{
       \setbox\mybox=\hbox spread 1.4em{\hfil$\scriptstyle#1$\hfil}
       \vbox{\offinterlineskip\copy\mybox
             \hbox to\wd\mybox{\rightarrowfill}}}}
\def\larrover#1{\mathrel{
       \setbox\mybox=\hbox spread 1.4em{\hfil$\scriptstyle#1$\hfil}
       \vbox{\offinterlineskip\copy\mybox
             \hbox to\wd\mybox{\leftarrowfill}}}}

\def\ontoover#1{\mathrel{
       \setbox\mybox=\hbox spread 1.4em{\hfil$\scriptstyle#1$\hfil}
       \vbox{\offinterlineskip\copy\mybox
             \hbox to\wd\mybox{\rightarrowfill\hskip-2.8mm
                               $\rightarrow$}}}}
\def\leftontoover#1{\mathrel{
       \setbox\mybox=\hbox spread 1.4em{\hfil$\scriptstyle#1$\hfil}
       \vbox{\offinterlineskip\copy\mybox
             \hbox to\wd\mybox{$\leftarrow$\hskip-2.8mm
                               \leftarrowfill}}}}
\def\longto{\longrightarrow}
\def\into{\hookrightarrow}
\def\onto{\ontoover{\ }}
\def\longonto{\ontoover{\ }}
\def\isoto{\arrover{\sim}}

\def\longinto{\lhook\joinrel\longrightarrow}

\usepackage[curve,matrix,arrow,cmtip]{xy}
\NoComputerModernTips

\def\myxymessage{\def\messagetext
   {Here an xy-pic diagram was omitted to speed up compilation . . . }
   \immediate\write16{\messagetext}
   \hbox{\bf \messagetext}}
\def\filxymatrix#1{\myxymessage}
\def\filxyarray#1{\myxymessage}
\newdir^{ (}{{}*!/-3pt/\dir^{(}}
\newdir_{ (}{{}*!/-3pt/\dir_{(}}
\newdir^{ )}{{}*!/+3pt/\dir^{)}}
\newdir_{ )}{{}*!/+3pt/\dir_{)}}

\def\rscript#1{\hbox to 0pt{$\scriptstyle#1$\hss}}

\let\oldbullet\bullet
\def\bullet{{\mathchoice{\oldbullet}%
                        {\oldbullet}%
                        {\scriptscriptstyle\oldbullet}%
                        {\oldbullet}}}

\newcommand{\bS}{\mathop{\overline{S}}\nolimits}
\newcommand{\bX}{\mathop{\overline{X}}\nolimits}
\newcommand{\bbX}{\mathop{\overline{\overline{X}}}\nolimits}
\newcommand{\CHeffM}{\mathop{CHM^{eff}(k)}\nolimits}
\newcommand{\CHM}{\mathop{CHM(k)}\nolimits}
\newcommand{\cone}{\mathop{\rm Cone}\nolimits}

\newcommand{\DeffgM}{\mathop{DM^{eff}_{gm}(k)}\nolimits}
\newcommand{\DeffgQM}{\mathop{DM^{eff}_{gm}(\BQ)}\nolimits}
\newcommand{\DgM}{\mathop{DM_{gm}(k)}\nolimits}

\newcommand{\uC}{\mathop{\underline{C}}\nolimits}
\newcommand{\Mgm}{\mathop{M_{gm}}\nolimits}
\newcommand{\Mcgm}{\mathop{M_{gm}^c}\nolimits}
\newcommand{\dMgm}{\mathop{\partial M_{gm}}\nolimits}

\newcommand{\jp}{\mathop{\widetilde{j}} \nolimits}
\newcommand{\Xp}{\mathop{\widetilde{X}} \nolimits}
\newcommand{\Yp}{\mathop{\widetilde{Y}} \nolimits}

\begin{document}

%%%%%%%%%%%%%%%%%%%%%%%%%%%%%%%%%%%%%%%%%%%%%%%%%%%%%%%%%%%%%%%%%%%%%%%
%
%  formatting

\hfuzz=3pt
\overfullrule=10pt                   % erzeugt schwarze Fehlerbalken

% The displayskip values were changed because \LaTeX does not react
% correctly to a \leqno: it should then use big skips, but doesn't.

\setlength{\abovedisplayskip}{6.0pt plus 3.0pt}
                               % preset 10.0pt plus 2.0pt minus 5.0pt
\setlength{\belowdisplayskip}{6.0pt plus 3.0pt}
                               % preset 10.0pt plus 2.0pt minus 5.0pt
\setlength{\abovedisplayshortskip}{6.0pt plus 3.0pt}
                               % preset 0.0pt plus 3.0pt
\setlength{\belowdisplayshortskip}{6.0pt plus 3.0pt}
                               % preset 6.0pt plus 3.0pt minus 3.0pt

\setlength{\baselineskip}{13.0pt}
                               % preset 12.0pt
\setlength{\lineskip}{0.0pt}
                               % preset 1.0pt
\setlength{\lineskiplimit}{0.0pt}
                               % preset 0.0pt

%%%%%%%%%%%%%%%%%%%%%%%%%%%%%%%%%%%%%%%%%%%%%%%%%%%%%%%%%%%%%%%%%%%%%%%
%
%  Title Page
%
%%%%%%%%%%%%%%%%%%%%%%%%%%%%%%%%%%%%%%%%%%%%%%%%%%%%%%%%%%%%%%%%%%%%%%%

\title{Chow motives without projectivity
\forget{
\footnotemark
\footnotetext{To appear in ....}
}
}
\author{\footnotesize by\\ \\
\mbox{\hskip-2cm
\begin{minipage}{6cm} \begin{center} \begin{tabular}{c}
J\"org Wildeshaus \footnote{
Partially supported by the \emph{Agence Nationale de la
Recherche}, project no.\ ANR-07-BLAN-0142 ``M\'ethodes \`a la
Voevodsky, motifs mixtes et G\'eom\'etrie d'Arakelov''. } \\[0.2cm]
\footnotesize LAGA\\[-3pt]
\footnotesize UMR 7539\\[-3pt]
\footnotesize Institut Galil\'ee\\[-3pt]
\footnotesize Universit\'e Paris 13\\[-3pt]
\footnotesize Avenue Jean-Baptiste Cl\'ement\\[-3pt]
\footnotesize F-93430 Villetaneuse\\[-3pt]
\footnotesize France\\
{\footnotesize \tt wildesh@math.univ-paris13.fr}
\end{tabular} \end{center} \end{minipage}
\hskip-2cm}
%\\[4cm]
%{\bf Preliminary version --- not for distribution}\\[1cm]
}
% In the final version we might want to fix the date:
\date{January 16, 2009}
\maketitle
%\quad \\[-1.7cm]
\begin{abstract}
In \cite{Bo},
Bondarko recently defined the notion of
weight structure,
and proved that the category $\DgM$ of geometrical motives
over a perfect field $k$,
as defined and studied by Voevodsky,
Suslin and Friedlander \cite{VSF},
is canonically equipped with such a structure. 
Building on this result, and under a condition
on the weights avoided by the boundary motive \cite{W1},
we describe a method to 
construct intrinsically in $\DgM$ a
motivic version of interior cohomology
of smooth, but possibly non-projective schemes. 
In a sequel to this work \cite{W3}, 
this method will be applied to Shimura varieties. \\


%\noindent
%{\bf R\'esum\'e~:} RESUME.\\
\end{abstract}

%\vfill

\bigskip
\bigskip
\bigskip

\noindent {\footnotesize Math.\ Subj.\ Class.\ (2000) numbers: 14F42
(14F20, 14F25, 14F30, 14G35, 18E30, 19E15).}

\eject
\tableofcontents

\bigskip
%\vspace*{0.5cm}

%\newpage
%\include{Intro}

%%%%%%%%%%%%%%%%%%%%%%%%%%%%%%%%%%%%%%%%%%%%%%%%%%%%%%%%%%%%%%%%%%%%%%%
%
%  Introduction
%
%%%%%%%%%%%%%%%%%%%%%%%%%%%%%%%%%%%%%%%%%%%%%%%%%%%%%%%%%%%%%%%%%%%%%%%

\setcounter{section}{-1}
\section{Introduction}
\label{Intro}

%%%%%%%%%%%%%%%%%%%%%%%%%%%%%%%%

%%%%%%%%%%%%%%%%%%%%%%%%%%%%%%%%

The full title of this work would be ``Approximation of motives
of varieties which are smooth, but not necessarily projective,
by Chow motives, using homological
rather than purely geometrical methods''.   
This might give a better idea of our program ---
but as a title, it appeared too long. So there. \\  

One way to place the problem historically
is to start with the question asked by Serre \cite[p.~341]{Se}, 
whether the ``virtual motive'' $\chi_c(X)$ of an arbitrary
variety $X$ over a fixed base field
$k$ can be (well) defined in the Grothendieck group
of the ca\-tegory $\CHeffM$ of 
effective \emph{Chow motives}. When $k$ admits resolution
of singula\-rities, Gillet
and Soul\'e \cite{GS} (see also \cite{GN} when $char(k) = 0$) 
provided an affirmative
answer. In fact, their solution
yields much more information:
they define the \emph{weight complex} $W(X)$ ($h_c(X)$ in \cite{GN})
in the category of complexes over $\CHeffM$,
well defined up to canonical homotopy equi\-valence.
By definition, $\chi_c(X)$ equals the class of
$W(X)$ in $K_0 \bigl( \CHeffM \bigr)$. 
Thus, given any representative 
\[
M_\bullet : \quad\quad\quad\quad
\ldots \longto M_n \longto M_{n-1} \longto \ldots
\]
of $W(X)$, we have the fomula $\chi_c(X) = \sum_n (-1)^n [M_n]$. \\

Consider the fully faithful embedding $\iota$ of $\CHeffM$
into $\DeffgM$, the category of effective \emph{geometrical motives},
as defined and studied by Voe\-vodsky, Suslin and Friedlander \cite{VSF}. 
As observed in \cite{GS},
localization for both $\chi_c(X)$ and the 
\emph{motive with compact support} $\Mcgm(X)$ of $X$ shows that 
the element $\chi_c(X)$ is mapped 
to the class of $\Mcgm(X)$ under $K_0(\iota)$. 
Gillet and Soul\'e went on and asked  \cite[p.~153]{GS} whether 
\[
K_0(\iota) : K_0 \bigl( \CHeffM \bigr) \longto K_0 \bigl( \DeffgM \bigr)
\]
is an isomorphism. \\

Just as Serre, Gillet and Soul\'e provoked much more than
just an affirmative answer to their question. 
Bondarko defined in \cite{Bo} the notion of \emph{weight structure} 
on a triangulated category $\CC$. He proved that
the inclusion of the \emph{heart} 
of the weight structure into $\CC$ induces
an isomorphism on the level of Grothendieck groups,
whenever the weight structure is bounded, and the heart
is pseudo-Abelian.
The definitions will be recalled in our Section~\ref{1};
let us just note that shift by $[m]$, for $m \in \BZ$, adds
$m$ to the weight of an object of $\CC$.  
According to one of the main results of \loccit \
(recalled in Theorem~\ref{1D}), $\DeffgM$ carries 
a canonical weight structure, which is indeed bounded,
and admits $\CHeffM$ as its heart. 
By definition, this means that among all geometrical motives,
Chow motives distinguish themselves as being the motives which are
pure of weight zero. \\

In particular, this gives an intrinsic characterization of objects
of the category $\DeffgM$ belonging to $\CHeffM$. We think of this insight as 
nothing less than revolutionary. \\ 

To come back to the beginning, the component $M_n$ of the weight complex
$W(X)$ can be considered as a ``$\Gr_n W(X)$'' with respect
to the weight structure. In the context of Bondarko's theory,
the formula for the ``virtual motive'' of $X$ thus reads
\[
\chi_c(X) = \sum_n (-1)^n \Gr_n W(X) \; .
\]
Basically, our approach to construct a Chow motive out of $X$, when
$X$ is smooth, is very simple: according to Bondarko (see Corollary~\ref{1E}),
the \emph{motive} $\Mgm(X)$ of $X$ is of weights $\le 0$. We would like
to consider $\Gr_0 \Mgm(X)$, the ``quotient'' of $\Mgm(X)$ of maximal weight
zero. \\   

However, one of the main subtleties of the notion of weight structure is that
``the'' \emph{weight filtration} of an object is almost never unique.
(In the example of the weight complex, this corresponds to the fact that
$W(X)$ is only well-defined up to homotopy.) 
Hence ``$\Gr_n \Mgm(X)$'' is not well-defined, and the above approach 
cannot work as stated. 
In fact, for any smooth compacti\-fication $\Xp$ of $X$, the Chow motive
$\Mgm(\Xp)$ occurs as ``$\Gr_0 \Mgm(X)$'' for a suitable weight filtration! \\

Our main contribution is to identify a criterion
assuring existence and unicity of a ``best choice'' of $\Gr_0 \Mgm(X)$. 
It is best understood in the abstract setting created by Bondarko, that is,
in the context of a weight structure on a triangulated category $\CC$.
Let $M$ be an object of $\CC$ of non-positive weights. Then Bondarko's 
axioms imply that if $M$ admits a weight filtration
in which the adjacent weight $-1$ does not occur, then the
associated $\Gr_0 M$ is unique up to unique isomorphism. 
In the motivic context, this observation 
leads to our criterion
on weights avoided by the \emph{boundary motive}
$\dMgm(X)$ of $X$, introduced and studied in \cite{W1}. 
The behavior of the realizations of the Chow motive $Gr_0 \Mgm(X)$
motivates its name: we chose to call it the \emph{interior motive} of $X$. \\
 
Let us now give a more detailed description of the content 
of this article.
Section~\ref{1} claims no originality whatsoever.
We give Bondarko's definition of weight structures
(Definition~\ref{1A}) and review the results from \cite{Bo}
needed in the sequel.
We treat particularly carefully 
the phenomenon which will turn
out to be the main theme of this article, namely 
\emph{the absence of certain weights}.
Thus, we introduce the 
central notion of \emph{weight filtration avoiding weights
$m,m+1,\ldots,n-1,n$}, for fixed integers $m \le n$. 
Bondarko's axioms imply that whenever such a
filtration exists, it 
behaves functorially (Proposition~\ref{2B}). In particular, 
for any fixed object, it is 
unique up to unique isomorphism (Corollary~\ref{2Ca}). 
We conclude Section~\ref{1} with Bondarko's application of his theory
to geometrical motives, which we already mentioned 
before (Theorem~\ref{1D}, Corollary~\ref{1E}). \\

Section~\ref{2} is the technical center of this article.
It is devoted to a further study of functoriality of
weight filtrations avoiding certain weights. 
We work in the context of an abstract weight
structure $w = (\CC_{w \le 0},\CC_{w \ge 0})$ on a triangulated category 
$\CC$. We first show (Proposition~\ref{2Cc})
that the inclusion $\iota_-$ of the heart 
$\CC_{w = 0}$ into the full sub-category $\CC_{w \le 0, \ne -1}$
of objects of non-positive weights $\ne -1$
admits a left adjoint 
\[
\Gr_0 : \CC_{w \le 0, \ne -1} \longto \CC_{w = 0} \; .
\]
There is a dual version of this statement for non-negative weights $\ne 1$.
We then consider the situation which will be of interest in our application
to motives. We thus fix a morphism 
$u: M_- \to M_+$ between objects 
$M_- \in \CC_{w \le 0}$ and $M_+ \in \CC_{w \ge 0}$, and
a cone $C[1]$ of $u$. While the axioms characterizing weight structures
easily show that $u$ can be factored through \emph{some} 
object of $\CC_{w = 0}$,
our aim is to do so \emph{in a canonical way}. We therefore formulate
Assumption~\ref{2D}: the object $C$ is without weights $-1$ and $0$.
Theorem~\ref{2main} states that
Assumption~\ref{2D} not only allows to factor $u$ as desired; in addition
this factorization is through an object which is simultaneously identified
with $\Gr_0 M_-$ and with $\Gr_0 M_+$. As a formal consequence 
of this, and of the functoriality properties of $\Gr_0$ 
from Proposition~\ref{2Cc},
we get a statement on abstract
factorization of $u$ (Corollary~\ref{2E}), whose rigidity may appear
surprising at first sight, 
given that we work in a triangulated category: whenever
$u: M_- \to N \to M_+$ factors through an object
$N$ of $\CC_{w = 0}$, then $\Gr_0 M_- = \Gr_0 M_+$ is canonically
identified with
a direct factor of $N$, admitting in addition a canonical direct complement. \\

The reader willing to turn directly to the 
application of these
results to geo\-metrical motives
may choose to skip Section~\ref{4}, in which 
Scholl's construction of motives for modular forms
\cite{Scho} is discussed at length. 
As in \loccit , we consider a self-product $X_n^r$ of the universal elliptic
curve $X_n$ over a modular curve. We 
show (Theorem~\ref{4Ca},
Corollary~\ref{4C}, Corollary~\ref{4Da})
that certain direct factors $\Mgm(X_n^r)^e$
of the motive $\Mgm(X_n^r)$ and $\Mcgm(X_n^r)^e$ of 
the motive with compact support $\Mcgm(X_n^r)$,
together with the canonical morphism $u: \Mgm(X_n^r)^e \to \Mcgm(X_n^r)^e$,
satisfy the conclusions of Theorem~\ref{2main}. In particular,
the Chow motives $\Gr_0 \Mgm(X_n^r)^e$ and $\Gr_0 \Mcgm(X_n^r)^e$
are defined, and
canonically isomorphic. In fact, they are both canonically
isomorphic to the motive denoted ${}^r_n \CW$ in \cite{Scho}. We insist on
giving a proof of the above statements which is independent of
the theory developed in Section~\ref{2}, that is, we prove the
conclusions of Theorem~\ref{2main} without first
checking Assumption~\ref{2D}. Instead, we use  
the detailed analysis from \cite{Scho} of the geometry of the boundary
of a ``good choice'' of smooth compactification of $X_n^r$. 
If one forgets about the language of weight structures, whose use could
not be completely avoided, Section~\ref{4} is thus technically independent
of the material preceding it. The reader may find an interest in
the re-interpretation of Scholl's construction in the context of 
Voevodsky's
geo\-metrical motives, which were not yet defined
at the time when \cite{Scho} was written. This concerns in particular
the exact triangle
\[
\Mgm \bigl( S_n^\infty \bigr) (r+1)[r+1] 
\longto \Mgm(X_n^r)^e 
\longto {}^r_n \CW
\longto \Mgm \bigl( S_n^\infty \bigr) (r+1)[r+2] 
\]  
from Corollary~\ref{4C} ($S_n^\infty :=$ the cuspidal locus of the modular
curve). As a by-product of Scholl's analysis, we prove that the triangle is 
defined using $\BZ[1 / (2n \cdot r!)]$-coefficients.
An interpretation of Beilinson's \emph{Eisenstein symbol} \cite{B}
in this context (namely, as a splitting of this triangle)
is clearly desirable. \\

Section~\ref{3} is devoted to 
the application of the
results from Section~\ref{2} to geo\-metrical motives. 
As $u: M_- \to M_+$, we take the canonical morphism
$\Mgm(X)^e \to \Mcgm(X)^e$, for
a fixed smooth variety $X$ over $k$, and a fixed idempotent $e$.
Then the role of the object $C$ from Section~\ref{2} is canonically
played by $\dMgm(X)^e$, the $e$-part of the boundary motive.
In this context,
Assumption~\ref{2D} reads as follows: the object $\dMgm(X)^e$ 
is without weights $-1$ and $0$ (Assumption~\ref{3A}).
Our main result Theorem~\ref{Main} is then simply the translation of
the sum of the results from Section~\ref{2} into this 
particular motivic context.
Thus, the Chow motive $\Gr_0 \Mgm(X)^e = \Gr_0 \Mcgm(X)^e$ 
is defined. Furthermore (Corollary~\ref{3B}),
we get the motivic version of abstract
factorization: whenever
$\Xp$ is a smooth compactification of $X$, then
$\Gr_0 \Mgm(X)^e$ is canoni\-cally
a direct factor of the Chow motive
$\Mgm(\Xp)$, with a canonical direct complement.
We then study the implications of these results for the 
Hodge theoretic and $\ell$-adic realizations. Theorems~\ref{3C}
and \ref{3D} state 
that they are equal to the respective $e$-parts of interior cohomology of $X$.
Abstract factorization allows to say more about the quality of the
Galois representation on the $\ell$-adic realization 
of $\Gr_0 \Mgm(X)^e$ (Theorem~\ref{3F}).
For example, simple semi-stable reduction of 
\emph{some} smooth compactification of $X$
implies that the representation is semi-stable.   \\

To conclude, we get back to Scholl's construction.
We show (Remark~\ref{3H}) that essentially 
all results of Section~\ref{4} can be deduced (just)
from Assumption~\ref{3A}, and the theory
developed in Section~\ref{3}. We consider that
in spite of the technical independence of Section~\ref{4}, 
there is a good reason to include that material
in this article: for higher dimensional Shimura
varieties, ``good choices'' of smooth compactifications 
as the one used in \cite{Scho} may simply not be available.
Therefore, the purely geometrical strategy of proof of
the results of Section~\ref{4} 
can safely be expected not to 
be generalizable. We think that a promising way to generalize is 
via a verification of Assumption~\ref{3A} by other than
purely geometrical means. We refer to \cite{W3} for the development
of such an alternative
in a context including that of (powers of universal 
elliptic curves over) modular curves. \\
 
Part of this work was done while I was enjoying a 
\emph{modulation de service pour les porteurs de projets de recherche},
granted by the \emph{Universit{\'e} Paris~13}. 
I wish to thank M.V.~Bondarko, F.~D\'eglise,
A.~Deitmar, O.~Gabber, B.~Kahn, A.~Mo\-krane,
C.~Soul\'e and J.~Tilouine for 
useful discussions, and the referee(s) for helpful suggestions. \\

{\bf Notation and conventions}: $k$ denotes a fixed perfect
base field, $Sch/k$  
the category of separated schemes of finite type over $k$, and  
$Sm/k \subset Sch/k$
the full sub-category of objects which are smooth over $k$.
When we assume $k$ to admit resolution of singularities,
then it will be 
in the sense of \cite[Def.~3.4]{FV}:
(i)~for any $X \in Sch/k$, there exists an abstract blow-up $Y \to X$
\cite[Def.~3.1]{FV} whose source $Y$ is in $Sm/k$,
(ii)~for any $X, Y \in Sm/k$, and any abstract blow-up $q : Y \to X$,
there exists a sequence of blow-ups 
$p: X_n \to \ldots \to X_1 = X$ with smooth centers,
such that $p$ factors through $q$. 
Let us note that the main reason for us to suppose $k$ to admit resolution
of singularities is to have the motive 
with compact support satisfy localization \cite[Prop.~4.1.5]{V}. \\   

As far as motives are concerned,
the notation of this paper follows
that of \cite{V}. We refer to \cite[Sect.~1]{W1} for a 
review of this notation, and in particular,
of the definition
of the categories $\DeffgM$ and $\DgM$ of (effective) geometrical
motives over $k$,
and of the motive $\Mgm(X)$ and the motive with compact support $\Mcgm(X)$
of $X \in Sch / k$. Let $F$ be a commutative flat $\BZ$-algebra,
i.e., a commutative unitary ring whose additive group is without torsion.
The notation $\DeffgM_F$ and $\DgM_F$ stands for the 
$F$-linear analogues of $\DeffgM$ and $\DgM$
defined in \cite[Sect.~16.2.4
and Sect.~17.1.3]{A}. 
Simi\-larly,
let us denote by $\CHeffM$ and $\CHM$ the categories opposite to the categories
of (effective) Chow motives, and by
$\CHeffM_F$ and $\CHM_F$ the pseudo-Abelian
completion of the category
$\CHeffM \otimes_\BZ F$ and $\CHM \otimes_\BZ F$,
respectively. Using \cite[Cor.~2]{V2} (\cite[Cor.~4.2.6]{V} if $k$ 
admits resolution of singularities), we canonically identify 
$\CHeffM_F$ and $\CHM_F$ with
a full additive sub-category of $\DeffgM_F$ and $\DgM_F$, respectively.

%%% Local Variables:
%%% mode: latex
%%% TeX-master: "head"
%%% End:

\bigskip
%\include{Sec1}

%%%%%%%%%%%%%%%%%%%%%%%%%%%%%%%%%%%%%%%%%%%%%%%%%%%%%%%%%%%%%%%%%%%%%%%
%
%  Section 1
%
%%%%%%%%%%%%%%%%%%%%%%%%%%%%%%%%%%%%%%%%%%%%%%%%%%%%%%%%%%%%%%%%%%%%%%%

\section{Weight structures}
\label{1}

%%%%%%%%%%%%%%%%%%%%%%%%%%%%%%%%

%%%%%%%%%%%%%%%%%%%%%%%%%%%%%%%%

In this section, we review definitions and
results of Bondarko's recent paper \cite{Bo}.

\begin{Def} \label{1A}
Let $\CC$ be a triangulated category. A \emph{weight structure on $\CC$}
is a pair $w = (\CC_{w \le 0} , \CC_{w \ge 0})$ of full 
sub-categories of $\CC$, such that, putting
\[
\CC_{w \le n} := \CC_{w \le 0}[n] \quad , \quad
\CC_{w \ge n} := \CC_{w \ge 0}[n] \quad \forall \; n \in \BZ \; ,
\]
the following conditions are satisfied.
\begin{enumerate}
\item[(1)] The categories
$\CC_{w \le 0}$ and $\CC_{w \ge 0}$ are 
Karoubi-closed: for any object $M$ of $\CC_{w \le 0}$ or
$\CC_{w \ge 0}$, any direct summand of $M$ formed in $\CC$
is an object of $\CC_{w \le 0}$ or
$\CC_{w \ge 0}$, respectively.
\item[(2)] (Semi-invariance with respect to shifts.)
We have the inclusions
\[
\CC_{w \le 0} \subset \CC_{w \le 1} \quad , \quad
\CC_{w \ge 0} \supset \CC_{w \ge 1}
\]
of full sub-categories of $\CC$.
\item[(3)] (Orthogonality.)
For any pair of objects $M \in \CC_{w \le 0}$ and $N \in \CC_{w \ge 1}$,
we have
\[
\Hom_{\CC}(M,N) = 0 \; .
\]
\item[(4)] (Weight filtration.)
For any object $M \in \CC$, there exists an exact triangle
\[
A \longto M \longto B \longto A[1]
\]
in $\CC$, such that $A \in \CC_{w \le 0}$ and $B \in \CC_{w \ge 1}$.
\end{enumerate}
\end{Def}

By condition~\ref{1A}~(2),
\[
\CC_{w \le n} \subset \CC_{w \le 0}
\]
for negative $n$, and
\[
\CC_{w \ge n} \subset \CC_{w \ge 0}
\]
for positive $n$. There are obvious analogues of the other conditions
for all the categories $\CC_{w \le n}$ and $\CC_{w \ge n}$. In particular,
they are all Karoubi-closed, and any object $M \in \CC$
is part of an exact triangle
\[
A \longto M \longto B \longto A[1]
\]
in $\CC$, such that $A \in \CC_{w \le n}$ and $B \in \CC_{w \ge n+1}$.
By a slight generalization of the terminology introduced in  
condition~\ref{1A}~(4), we shall refer to any such exact triangle
as a weight filtration of $M$.

\begin{Rem} \label{1Aa}
(a)~Our convention concerning the sign of the weight is actually opposite 
to the one from \cite[Def.~1.1.1]{Bo}, i.e., we exchanged the
roles of $\CC_{w \le 0}$ and $\CC_{w \ge 0}$. \\[0.1cm]
(b)~Note that in condition~\ref{1A}~(4), ``the'' weight filtration is not
assumed to be unique. \\[0.1cm]
(c)~Recall the notion of \emph{$t$-structure} on a triangulated category
$\CC$ \cite[D\'ef.~1.3.1]{BBD}. 
It consists of a pair $t = (\CC^{t \le 0},\CC^{t \ge 0})$ of full
sub-categories sa\-tisfying formal analogues of conditions~\ref{1A}~(2)--(4),
but putting
\[
\CC^{t \le n} := \CC^{t \le 0}[-n] \quad , \quad
\CC^{t \ge n} := \CC^{t \ge 0}[-n] \quad \forall \; n \in \BZ \; .
\]
Note that in the context of $t$-structures,
the analogues of the exact triangles in \ref{1A}~(4) are then
unique up to unique isomorphism, and that the analogue of 
condition~\ref{1A}~(1) is formally implied by the others.
\end{Rem}

The following is contained in \cite[Def.~1.2.1]{Bo}.
 
\begin{Def} \label{1B}
Let $w = (\CC_{w \le 0} , \CC_{w \ge 0})$ be a weight structure on $\CC$.
The \emph{heart of $w$} is the full additive sub-category $\CC_{w = 0}$
of $\CC$ whose objects lie 
both in $\CC_{w \le 0}$ and in $\CC_{w \ge 0}$.
\end{Def}

Among the basic properties developed in \cite{Bo}, 
let us note the following.

\begin{Prop} \label{1C}
Let $w = (\CC_{w \le 0} , \CC_{w \ge 0})$ be a weight structure on $\CC$,
\[
L \longto M \longto N \longto L[1]
\]
an exact triangle in $\CC$. \\[0.1cm]
(a)~If both $L$ and $N$ belong to $\CC_{w \le 0}$, then so does $M$. \\[0.1cm]
(b)~If both $L$ and $N$ belong to $\CC_{w \ge 0}$, then so does $M$. 
\end{Prop}

\begin{Proof}
This is the content of \cite[Prop.~1.3.3~3]{Bo}.
\end{Proof}

The reader may wonder whether there is an easy criterion on
a given sub-category of a
triangulated category to be the heart of a suitable weight structure.
Bondarko has results \cite[Thm.~4.3.2]{Bo}
answering this question. For our purposes, the result with the
most restrictive finiteness condition will be sufficient.

\begin{Prop} \label{1Ca}
Let $\CH$ be a full additive sub-category of a triangulated category $\CC$.
Suppose that $\CH$ generates $\CC$, i.e., $\CC$ is the smallest full
triangulated sub-category containing $\CH$. \\[0.1cm]
(a)~If there is a weight structure on $\CC$ whose heart contains $\CH$,
then it is unique. In this case, the heart is equal to the \emph{Karoubi
envelope} of $\CH$, i.e., the category of retracts of $\CH$ in $\CC$. \\[0.1cm]
(b)~The following conditions are equivalent.
\begin{enumerate}
\item[(i)] There is a weight structure on $\CC$ whose heart
contains $\CH$.
\item[(ii)] $\CH$ is \emph{negative}, i.e.,
\[
\Hom_{\CC} \bigl( A , B [i] \bigr) = 0
\]
for any two objects $A$, $B$ of $\CH$, and any integer $i > 0$.
\end{enumerate}
\end{Prop}

\begin{Proof}
Condition (ii) on $\CH$ is clearly necessary for $\CH$ to belong to the heart,
given orthogonality~\ref{1A}~(3).
As for (a), note that by Proposition~\ref{1C},
there is only one possible definition of
the category $\CC_{w \le 0}$ (resp.\ $\CC_{w \ge 0}$): 
it is necessarily the full sub-category of
successive extensions of objects of the form $A[n]$, for $A \in \CH$
and $n \le 0$ (resp.\ $n \ge 0$).

The main point is to show
that under condition (ii), the above
construction indeed yields a weight structure on $\CC$.   
We refer to \cite[Thm.~4.3.2~II]{Bo} for details.
\end{Proof}

For the rest of this section,
we consider a fixed weight structure $w$ on
a triangulated category $\CC$.

\begin{Def} \label{2A}
Let $M \in \CC$, and $m \le n$ two integers (which may be identical). 
A \emph{weight filtration of $M$ avoiding weights $m,m+1,\ldots,n-1,n$}
is an exact triangle
\[ 
M_{\le {m-1}} \longto M \longto M_{\ge {n+1}} \longto M_{\le m-1}[1]
\]
in $\CC$, with $M_{\le {m-1}} \in \CC_{w \le m-1}$
and $M_{\ge {n+1}} \in \CC_{w \ge {n+1}} \ $. \\[0.1cm]
\end{Def}

The following observation is vital.

\begin{Prop} \label{2B}
Assume that $m \le n$, and
that $M, N \in \CC$ admit weight filtrations 
\[
M_{\le {m-1}} \stackrel{x_-}{\longto} M 
\stackrel{x_+}{\longto} M_{\ge {n+1}} \longto M_{\le m-1}[1]
\]
and
\[
N_{\le {m-1}} \stackrel{y_-}{\longto} N 
\stackrel{y_+}{\longto} N_{\ge {n+1}} \longto N_{\le m-1}[1]
\]
avoiding weights $m,\ldots,n$. Then any morphism $M \to N$ in $\CC$
extends uniquely to a morphism of exact triangles 
\[
\vcenter{\xymatrix@R-10pt{
        M_{\le {m-1}} \ar[r] \ar[d] &
        M \ar[r] \ar[d] &
        M_{\ge {n+1}} \ar[r] \ar[d] &
        M_{\le m-1}[1] \ar[d] \\
        N_{\le {m-1}} \ar[r] &
        N \ar[r] &
        N_{\ge {n+1}} \ar[r] &
        N_{\le m-1}[1] 
\\}}
\]
\end{Prop}

\begin{Proof}
This follows from \cite[Lemma~1.5.1~2]{Bo}.
For the convenience of the reader, let us recall the proof. 
Let $\alpha \in \Hom_\CC(M,N)$.
The composition $y_+ \circ \alpha \circ x_- : M_{\le {m-1}} \to N_{\ge {n+1}}$  
is zero by orthogonality~\ref{1A}~(3): $m-1$ is strictly smaller than $n+1$. 
Hence $\alpha \circ x_-$ factors through $N_{\le {m-1}}$.
We claim that this factorization is unique. Indeed, the 
error term comes from 
$\Hom_\CC( M_{\le m-1} , N_{\ge {n+1}}[-1] )$. But this group is trivial, 
thanks to orthogonality, and
our assumption on the weights: the object $N_{\ge {n+1}}[-1]$ lies in
\[
\CC_{w \ge n+1}[-1] = \CC_{w \ge n} \; ,
\]
and $m-1$ is still stricly smaller than $n$.
Similarly, the composition $y_+ \circ \alpha$ factors uniquely through 
$M_{\ge {n+1}}$.
\end{Proof}

\begin{Rem} \label{2C}
Note that the hypothesis of Proposition~\ref{2B}
does not imply unicity of weight filtrations
in the (more general) sense of \ref{1A}~(4). 
For example, assume that $m = n = -1$, and let
\[
(\ast) \quad\quad
M_{\le -2} \longto M \longto M_{\ge 0} \longto M_{\le -2}[1]
\]
be a weight filtration avoiding weight $-1$.
Choose any object $M_0$ in $\CC_{w = 0}$ and
replace $M_{\ge 0}$ by 
$M_0 \oplus M_{\ge 0}$, and $M_{\le -2}$ by $M_0[-1] \oplus M_{\le -2}$.
Arguing as in the proof of Proposition~\ref{2B}, 
one shows that any weight filtration of $M$ 
is isomorphic to one obtained 
in this way. Thus, the exact triangle $(\ast)$ satisfies
a minimality property among all weight filtrations of $M$. 
\end{Rem}

\begin{Cor} \label{2Ca}
Assume that $m \le n$. Then if
$M \in \CC$ admits a weight filtration 
avoiding weights $m,\ldots,n$, it is unique up to unique isomorphism. 
\end{Cor}

\begin{Def}
Assume that $m \le n$. We say that
$M \in \CC$ \emph{does not have weights $m,\ldots,n$},
or that $M$ \emph{is without weights $m,\ldots,n$},
if it admits a weight filtration 
avoiding weights $m,\ldots,n$.
\end{Def}

Let us now state what we consider as one of the
main results of \cite{Bo}.

\begin{Thm} \label{1Cb}
Assume that the triangulated category $\CC$ is generated by its
heart $\CC_{w = 0}$. \\[0.1cm]
(a)~The pseudo-Abelian completion $\CC_{w = 0}'$ of $\CC_{w = 0}$
generates the pseudo-Abelian completion $\CC'$ of $\CC$. \\[0.1cm]
(b)~There is a weight structure $w'$ on $\CC'$, uniquely characterized
by any of the following conditions.  
\begin{enumerate}
\item[(i)] The weight structure $w'$ extends $w$.
\item[(ii)] The heart of $w'$ equals $\CC_{w = 0}'$.
\item[(iii)] The heart of $w'$ contains $\CC_{w = 0}'$.
\end{enumerate}
\end{Thm}

\begin{Proof}
This is \cite[Prop.~5.2.2]{Bo}. Let us describe the main steps
of the proof.
Recall that by \cite[Thm.~1.5]{BaSc}, the category $\CC'$ is indeed
triangulated. The criterion from Proposition~\ref{1Ca}
implies the existence of a weight structure $w'$
on the full triangulated sub-category $\CD$ of $\CC'$ generated
by $\CC_{w = 0}'$ (hence containing $\CC$), 
and uniquely characterized by condition (iii),
hence also by (i) or (ii).
The claim then follows from \cite[Lemma~5.2.1]{Bo}, which states that
$\CD$ is pseudo-Abelian, and hence equal to $\CC'$.
\end{Proof}

\begin{Rem}
Note that given Proposition~\ref{1Ca}, part (b) of Theorem~\ref{1Cb}
follows formally from its part (a). One may see Theorem~\ref{1Cb}~(a)
as a gene\-ralization of \cite[Cor.~2.12]{BaSc}, which states that the pseudo-Abelian
completion of the bounded derived 
category $D^b(\CA)$ of an exact 
category $\CA$ equals the bounded derived category $D^b(\CA')$ of the
pseudo-Abelian completion $\CA'$ of $\CA$.
\end{Rem}

For our purposes, the main application of the preceding
is the following (cmp.\ \cite[Sect.~6]{Bo}).

\begin{Thm} \label{1D}
Let $F$ be a commutative flat $\BZ$-algebra,
and assume $k$ to admit resolution of singularities. \\[0.1cm]
(a)~There is a canonical weight structure on the category $\DeffgM_F$.
It is uniquely characterized by the requirement that its heart equal 
$\CHeffM_F$. \\[0.1cm]
(b)~There is a canonical weight structure on the category $\DgM_F$,
extending the weight structure from (a).
It is uniquely characterized by the requirement that its heart equal 
$\CHM_F$. \\[0.1cm]
(c)~Statements (a) and (b) hold without assuming resolution of singularities  
provided $F$ is a $\BQ$-algebra.
\end{Thm}

\begin{Proof}
For $F = \BZ$ and $k$ of characteristic zero,
this is the content of \cite[Sect.~6.5 and 6.6]{Bo}:

(1)~As in \cite{Bo}, denote by $DM^s$ the full triangulated sub-category
of $\DeffgM$
generated by the motives $\Mgm(X)$ \cite[Def.~2.1.1]{V} 
of objects $X$ of $Sm/k$, by
$J_0$ the full additive sub-category of $DM^s$ generated by $\Mgm(X)$ for
$X$ smooth and projective, and by $J_0'$ the Karoubi envelope
of $J_0$. Thus, $\DeffgM$ is the pseudo-Abelian 
completion of $DM^s$, and $\CHeffM$ is the pseudo-Abelian completion of
both $J_0$ and $J_0'$. 

(2)~We need two of the main results from \cite{V}.
First, by \loccit , Cor.~3.5.5, the additive category $J_0$
generates the triangulated category $DM^s$.
Next, by \loccit , Cor.~4.2.6, the category $J_0$ is negative:
\[
\Hom_{DM^s} \bigl( A , B [i] \bigr) = 0
\]
for any two objects $A$, $B$ of $J_0$, and any integer $i > 0$.

(3)~By Proposition~\ref{1Ca}, 
there is a weight structure on $DM^s$, uniquely characterized by
the fact that $J_0$ is contained in the heart. Furthermore, the heart
equals the Karoubi envelope $J_0'$. 

(4)~By Theorem~\ref{1Cb}, the pseudo-Abelian
completion $\CHeffM$ of $J_0$ generates the pseudo-Abelian completion
$\DeffgM$ of $DM^s$ (let us remark that this is stated,
but not proved in \cite[Cor.~3.5.5]{V}). Thus, 
part (a) of our claim holds for $F = \BZ$.

(5)~Recall that $\CHM$ and $\DgM$ are obtained from $\CHeffM$ and $\DeffgM$
by inverting an object, namely the Tate object $T$, 
with respect to the tensor structures. 
Hence $\CHM$ ge\-ne\-rates the triangulated category $\DgM$.
Its negativity follows formally from that of $\CHeffM$: indeed,
for two objects $A$ and $B$ of $\CHM$, and any integer $i > 0$,
the group $\Hom_{\DgM} (A,B[i])$ is by definition the direct limit over 
large integers $r$ of the groups
\[
\Hom_{\DeffgM} \bigl( A \otimes T^{\otimes r},
                          B \otimes T^{\otimes r}[i] \bigr) \; ,
\]
which are all zero by part~(a). Thus,
we may again apply Proposition~\ref{1Ca}. The resulting weight structure
extends the one on $\DeffgM$: in fact, its restriction to 
$\DeffgM$ is a weight structure, whose heart equals $\CHeffM$. 
This proves part (b) of our claim for $F = \BZ$. 

If $F$ is flat over $\BZ$, then the same proof works.
(1')~Replace $DM^s$ by the full $F$-linear triangulated
sub-category $DM^s_F$ of $\DeffgM_F$
generated by the motives $\Mgm(X)$ of objects $X$ of $Sm/k$,
and $J_0$ by $J_0 \otimes_\BZ F$.
(2')~The two results cited in (2) formally imply that $J_0 \otimes_\BZ F$
generates $DM^s_F$, and that $J_0 \otimes_\BZ F$ is negative.
Steps (3') and (4') are formally identical to (3) and (4), proving part (a)
of the claim. Step (5') shows part (b), once we observe that 
$\CHM_F$ and $\DgM_F$ are obtained from $\CHeffM_F$ and $\DeffgM_F$
by inverting the Tate object.

As for part (c) of our claim, everything reduces to showing analogues of
the two statements made in step (2). 
By \cite[Cor.~18.1.1.2]{A},
the additive ca\-te\-gory $J_0 \otimes_\BZ F$
generates the triangulated category $DM^s \otimes_\BZ F$. The argument
uses alterations \`a la de Jong; since this involves finite extensions
of fields, whose degrees need to be inverted, one requires
$F$ to be a $\BQ$-algebra.
The generalization of \cite[Cor.~4.2.6]{V} to arbitrary fields
\cite[Cor.~2]{V2} shows that the category $J_0$ is negative.
Hence so is $J_0 \otimes_\BZ F$.
\end{Proof}

The following is the content of \cite[Thm.~6.2.1~1 and 2]{Bo1}.

\begin{Cor} \label{1E}
Assume $k$ to admit resolution of singularities.
Let $X$ in $Sch/k$ be of (Krull) dimension $d$. \\[0.1cm]
(a)~The motive with compact support $\Mcgm(X)$ lies in 
\[
\DeffgM_{w \ge 0} \cap \DeffgM_{w \le d} \; . 
\]
(b)~If $X \in Sm/k$,
then the motive $\Mgm(X)$ 
lies in 
\[
\DeffgM_{w \ge -d} \cap \DeffgM_{w \le 0} \; .
\]
\end{Cor}

\begin{Proof}
(a)~We proceed by induction on $d$. If
$d = 0$, then $\Mcgm(X)$ is an effective
Chow motive, hence of weight $0$ by Theorem~\ref{1D}~(a).

For $d \ge 1$, Nagata's theorem on the existence of a compactification 
of $X$, and resolution of singularities imply that there is an open
dense sub-scheme $U$ of $X$ admitting a smooth compactification $\Xp$.
Denote by $Z$ the complement of $U$ in $X$, and by $Y$ the complement 
of $U$ in $\Xp$ (both with the reduced scheme structure).
By localization for the motive with compact support \cite[Prop.~4.1.5]{V},
there are exact triangles
\[
\Mcgm(Z) \longto \Mcgm(X) \longto \Mcgm(U) \longto \Mcgm(Z)[1] \; .
\]
and 
\[
\Mcgm(Y) \longto \Mcgm(\Xp) \longto \Mcgm(U) \longto \Mcgm(Y)[1] \; .
\]
By induction, 
\[
\Mcgm(Y), \Mcgm(Z) \in \DeffgM_{w \ge 0} \cap \DeffgM_{w \le d-1} \; .
\]
Therefore, 
\[
\Mcgm(Y)[1], \Mcgm(Z) \in \DeffgM_{w \ge 0} \cap \DeffgM_{w \le d} \; .
\]
Given that $\Mcgm(\Xp)$ is of weight $0$,
Proposition~\ref{1C} shows first that 
\[
\Mcgm(U) \in \DeffgM_{w \ge 0} \cap \DeffgM_{w \le d} \; ,
\]
and then that
\[
\Mcgm(X) \in \DeffgM_{w \ge 0} \cap \DeffgM_{w \le d} \; .
\]
(b)~By \cite[Thm.~4.3.7~1 and 2]{V}, the category $\DgM$ is
rigid tensor triangulated. The claim thus follows formally from (a),
and from the following observations:
(i)~assuming (as we may) $X$ to be of pure dimension $d$,
the motive $\Mgm(X)$ is dual to $\Mcgm(X)(-d)[-2d]$
\cite[Thm.~4.3.7~3]{V},
(ii)~the object $\BZ(-d)[-2d]$ is a Chow motive,
(iii)~the heart of the weight structure on 
$\DgM$ is stable under
duality, hence for any natural number $n$, induction on $n$ shows that
the dual of an object of the intersection
$\DgM_{w \ge 0} \cap \DgM_{w \le n}$ belongs to 
$\DgM_{w \ge -n} \cap \DgM_{w \le 0}$,
(iv)~the weight structure on $\DeffgM$ is induced from the
weight structure on $\DgM$ (Theorem~\ref{1D}~(b)).
\end{Proof}

\begin{Rem}
(a)~Corollary~\ref{1E}~(a) and its proof should be compared to
the construction of the weight complex $W(X)$ 
from \cite[Sect.~2.1]{GS}. Let $\jp: \Yp_\bullet \to \Xp_\bullet$
be a smooth hyper-envelope (in the sense of \loccit )
of a closed immersion $Y \into \bX$ of proper schemes
whose complement equals $X$.
Both $\Yp_\bullet$ and $\Xp_\bullet$ give rise to complexes 
$\BZ \Yp_\bullet$ and $\BZ \Xp_\bullet$ in the category denoted
$\BZ {\bf V}$ in \loccit , i.e., the $\BZ$-linearized category
associated to the category ${\bf V}$ of smooth proper schemes
over $k$. Hence we may form the complex $\cone(\jp)$ 
in $\BZ {\bf V}$. Applying the functor $\uC \circ L$ 
\cite[pp.~207, 223--224]{V},
we get a complex of Nisnevich sheaves with transfers whose cohomology
sheaves are homotopy invariant. 
On the one hand,
it should be possible to employ 
\cite[Thm.~4.1.2]{V} in order to show that the complex $\uC (L(\cone(\jp)))$
represents the motive with compact support
$\Mcgm(X)$ in $\DeffgM$. On the other hand, by definition 
\cite[Sect.~2.1]{GS}, the (opposite of the) complex $\Mgm (\cone(\jp))$
represents $W(X)$ in the homotopy category over $\CHeffM$.
The statement  
\[
\Mcgm(X) \in \DeffgM_{w \ge 0} \cap \DeffgM_{w \le d} 
\]
from Corollary~\ref{1E}~(a) should be compared to 
\cite[Thm.~2~(i)]{GS}. \\[0.1cm]
(b)~Similarly, the construction from 
\cite[Thm.~(5.10)~(3)]{GN} of the object $h(X)$ 
should be compared to Corollary~\ref{1E}~(b).
\end{Rem}

\begin{Cor} \label{1F}
Assume $k$ to admit resolution of singularities.
Suppose given a direct factor $M$ of $\Mgm(X)$,
for $X \in Sm/k$, 
which is abstractly isomorphic to a direct factor of $\Mcgm(Y)$,
for some $Y \in Sch/k$. 
Then $M$ is an effective Chow motive.
\end{Cor}

%%% Local Variables:
%%% mode: latex
%%% TeX-master: "head"
%%% End:

\bigskip
%\include{Sec2}

%%%%%%%%%%%%%%%%%%%%%%%%%%%%%%%%%%%%%%%%%%%%%%%%%%%%%%%%%%%%%%%%%%%%%%%
%
%  Section 2
%
%%%%%%%%%%%%%%%%%%%%%%%%%%%%%%%%%%%%%%%%%%%%%%%%%%%%%%%%%%%%%%%%%%%%%%%

\section{Weight zero}
\label{2}

%%%%%%%%%%%%%%%%%%%%%%%%%%%%%%%%

%%%%%%%%%%%%%%%%%%%%%%%%%%%%%%%%

Throughout this section, we fix a weight structure 
$w = (\CC_{w \le 0} , \CC_{w \ge 0})$ on
a triangulated category $\CC$.
Anticipating the situation which will be of interest in 
our applications,
we formulate Assumption~\ref{2D} on the
cone of a morphism $u$ in $\CC$.
As we shall see (Theorem~\ref{2main}), this hypothesis
ensures in particular unique 
factorization of $u$ through an object of the heart $\CC_{w = 0}$. 

\begin{Def} \label{2Cb}
Denote by $\CC_{w \le 0, \ne -1}$ the full sub-category of $\CC_{w \le 0}$
of objects without weight $-1$, and
by $\CC_{w \ge 0, \ne 1}$ the full sub-category of $\CC_{w \ge 0}$
of objects without weight $1$. 
\end{Def}

\begin{Prop} \label{2Cc}
(a)~The inclusion of the heart 
$\iota_-: \CC_{w = 0} \into \CC_{w \le 0, \ne -1}$
admits a left adjoint 
\[
\Gr_0 : \CC_{w \le 0, \ne -1} \longto \CC_{w = 0} \; .
\]
On objects, it is given by sending $M$ to the term $M_{\ge 0}$
of a weight filtration
\[
M_{\le -2} \longto M \longto M_{\ge 0} \longto M_{\le -2}[1]
\]
avoiding weight $-1$. The composition $\Gr_0 \circ \iota_-$
equals the identity on $\CC_{w = 0}$. \\[0.1cm]
(b)~The inclusion of the heart 
$\iota_+: \CC_{w = 0} \into \CC_{w \ge 0, \ne 1}$
admits a right adjoint 
\[
\Gr_0 : \CC_{w \ge 0, \ne 1} \longto \CC_{w = 0} \; .
\]
On objects, it is given by sending $M$ to the term $M_{\le 0}$
of a weight filtration
\[
M_{\le 0} \longto M \longto M_{\ge 2} \longto M_{\le 0}[1]
\]
avoiding weight $1$. The composition $\Gr_0 \circ \iota_+$
equals the identity on $\CC_{w = 0}$.
\end{Prop}

\begin{Proof}
Given Proposition~\ref{2B}
and Corollary~\ref{2Ca}, all that remains to be proved is that
the objects $M_{\ge 0}$ (in (a)) resp.\ $M_{\le 0}$ (in (b))
actually do lie in $\CC_{w = 0}$. 
But this follows from Proposition~\ref{1C}.
\end{Proof}

Let us now fix the following data. 
\begin{enumerate}
\item[(1)] A morphism $u: M_- \to M_+$ in $\CC$ between 
$M_- \in \CC_{w \le 0}$ and $M_+ \in \CC_{w \ge 0} \ $.
\item[(2)] An exact triangle
\[
C \stackrel{v_-}{\longto} M_- \stackrel{u}{\longto}
M_+ \stackrel{v_+}{\longto} C[1] 
\]
in $\CC$. Thus, the object $C[1]$ is a fixed choice of cone of $u$. 
\end{enumerate}

We make the following rather restrictive hypothesis.

\begin{Ass} \label{2D}
The object $C$ is without weights $-1$ and $0$, i.e., it
admits a weight filtration 
\[ 
C_{\le -2} \stackrel{c_-}{\longto} C 
\stackrel{c_+}{\longto} C_{\ge 1} \stackrel{\delta_C}{\longto} C_{\le -2}[1]
\]
avoiding weights $-1$ and $0$.
\end{Ass}

The validity of this assumption is independent of the choice of $C$.
Here is our main technical tool.

\begin{Thm} \label{2main}
Fix the data (1), (2), and suppose Assumption~\ref{2D}. \\[0.1cm]
(a)~The object $M_-$ is without weight $-1$,
and $M_+$ is without weight $1$. \\[0.1cm]
(b)~The morphisms $v_- \circ c_-: C_{\le -2} \to M_-$ and
$\pi_0: M_- \to \Gr_0 M_-$ resp.\ $i_0: \Gr_0 M_+ \to M_+$
and $(c_+[1]) \circ v_+: M_+ \to C_{\ge 1}[1]$
can be canonically extended to exact triangles 
\[
(3) \quad\quad
C_{\le -2} \stackrel{v_- c_-}{\longto} M_- 
\stackrel{\pi_0}{\longto} \Gr_0 M_- \stackrel{\delta_-}{\longto} C_{\le -2}[1]
\]
and
\[
(4) \quad\quad
C_{\ge 1} \stackrel{\delta_+}{\longto} \Gr_0 M_+ \stackrel{i_0}{\longto} M_+ 
\stackrel{(c_+[1]) v_+}{\longto} C_{\ge 1}[1] \; .
\]
Thus, $(3)$ is a weight filtration of $M_-$ avoiding weight $-1$,
and $(4)$ is a weight filtration of $M_+$ avoiding weight $1$. \\[0.1cm]
(c)~There is a canonical isomorphism $\Gr_0 M_- \isoto \Gr_0 M_+$. 
As a morphism, it is uniquely determined by the property
of making the diagram
\[
\vcenter{\xymatrix@R-10pt{
        M_- \ar[r]^-{u} \ar[d]_{\pi_0} &
        M_+ \\
        \Gr_0 M_- \ar[r] &
        \Gr_0 M_+ \ar[u]_{i_0}
\\}}
\]
commute. Its inverse makes the diagram
\[
\vcenter{\xymatrix@R-10pt{
        C_{\ge 1} \ar[r]^-{\delta_C} \ar[d]_{\delta_+} &
        C_{\le -2}[1] \\
        \Gr_0 M_+ \ar[r] &
        \Gr_0 M_- \ar[u]_{\delta_-}
\\}}
\]
commute. 
\end{Thm} 

\begin{Proof}
We start by choosing and fixing exact triangles 
\[
(3') \quad\quad
C_{\le -2} \stackrel{v_- c_-}{\longto} M_- 
\stackrel{\pi_0'}{\longto} G_- \stackrel{\delta_-}{\longto} C_{\le -2}[1]
\]
and 
\[
(4') \quad\quad
C_{\ge 1} \stackrel{\delta_-}{\longto} G_+ \stackrel{i_0'}{\longto} M_+ 
\stackrel{(c_+[1]) v_+}{\longto} C_{\ge 1}[1] \; .
\]
Thus, $G_-$ is a cone of $v_- \circ c_-$,
and $G_+$ a cone of $c_+ \circ v_+[-1]$.
Observe first that by Proposition~\ref{1C}, 
\[
G_- \in \CC_{w \le 0} \quad \text{and} \quad 
G_+ \in \CC_{w \ge 0} \; .
\]
Given this, the existence of \emph{some} isomorphism
$G_- \cong G_+$ clearly implies parts (a) and (b) of our statement. 

Let us now show that there is an
isomorphism $\alpha: G_- \isoto G_+$
making the diagrams
\[
\vcenter{\xymatrix@R-10pt{
        M_- \ar[r]^-{u} \ar[d]_{\pi_0'} &
        M_+ \\
        G_- \ar[r]^{\alpha} &
        G_+ \ar[u]_{i_0'}
\\}}
\]
and
\[
\vcenter{\xymatrix@R-10pt{
        C_{\ge 1} \ar[r]^-{\delta_C} \ar[d]_{\delta_+} &
        C_{\le -2}[1] \\
        G_+ \ar[r]^{\alpha^{-1}} &
        G_- \ar[u]_{\delta_-}
\\}}
\]
commute. To see this, consider the following.
\[
\vcenter{\xymatrix@R-10pt{
        M_- \ar[dd]_{u}^{[1]} & &
        C_{\le -2} \ar[ll]_-{v_- \circ c_-} \ar[dl] \\
        & C \ar[ul] \ar[dr]^-{[1]} & \\
        M_+[-1] \ar[ur] \ar[rr]^{[1]}_-{c_+ \circ v_+[-1]} & &
        C_{\ge 1}[-1] \ar[uu]_{\delta_C[-1]}
\\}}
\]
The three arrows marked $[1]$ link the source to the shift by $[1]$
of the target,
the upper and lower triangles are commutative, and the left and right
triangles are exact. In the terminology of \cite[Sect.~1]{BBD},
this is a \emph{calotte inf\'erieure}, which thanks to the axiom TR4'
of triangulated categories in the formulation of \cite[1.1.6]{BBD}
can be completed to an octahedron. In particular, its contour
\[
\vcenter{\xymatrix@R-10pt{
        M_- \ar[dd]_{u}^{[1]} & &
        C_{\le -2} \ar[ll]_-{v_- \circ c_-} \\
        & & \\
        M_+[-1] \ar[rr]^{[1]}_-{c_+ \circ v_+[-1]} & &
        C_{\ge 1}[-1] \ar[uu]_{\delta_C[-1]}
\\}}
\]
is part of a \emph{calotte sup\'erieure}: 
\[
\vcenter{\xymatrix@R-10pt{
        M_- \ar[dd]_{u}^{[1]} \ar[dr]_-{[1]} & &
        C_{\le -2} \ar[ll]_-{v_- \circ c_-} \\
        & G \ar[ur] \ar[dl] & \\
        M_+[-1] \ar[rr]^{[1]}_-{c_+ \circ v_+[-1]} & &
        C_{\ge 1}[-1] \ar[uu]_{\delta_C[-1]} \ar[ul]
\\}}
\]
Here, the upper and lower triangles are exact, and the left and right
triangles are commutative. Hence the same object $G[1]$ can be chosen as
cone of $v_- \circ c_-$ and of $c_+ \circ v_+[-1]$, and in addition,
such that the morphisms $u$ and $\delta_C$ factor through $G[1]$.
For our fixed choices of cones, this means precisely 
that there is an isomorphism $G_- \cong G_+$
factorizing both $u$ and $\delta_C$.

As observed before, this implies (a) and (b).
It remains to show the unicity statement from (c). 
For this, use the exact triangle (3) and apply 
Proposition~\ref{2Cc}~(b) to see that
\[
\Hom_{\CC} \bigl( M_- , M_+ \bigr) 
= \Hom_{\CC} \bigl( \Gr_0 M_- , M_+ \bigr) 
= \Hom_{\CC_{w = 0}} \bigl( \Gr_0 M_- , \Gr_0 M_+ \bigr) \; . 
\]
Under this identification, a morphism $M_- \to M_+$ is sent
to its unique factorization $\Gr_0 M_- \to \Gr_0 M_+$.
\end{Proof}

Given Theorem~\ref{2main}~(c), we may and do identify 
$\Gr_0 M_-$ and $\Gr_0 M_+$.

\begin{Cor} \label{2E}
Fix the data (1), (2), and suppose Assumption~\ref{2D}. 
Let $M_- \to N \to M_+$ be a factorization of $u$ through an object
$N$ of $\CC_{w = 0}$. Then $\Gr_0 M_- = \Gr_0 M_+$ is canonically
identified with
a direct factor of $N$, admitting a canonical direct complement.
\end{Cor}

\begin{Proof}
By Proposition~\ref{2Cc}, the morphism $M_- \to N$ factors uniquely
through $\Gr_0 M_-$, and $N \to M_+$ factors uniquely through $\Gr_0 M_+$.
The composition $\Gr_0 M_- \to N \to \Gr_0 M_+$ is 
therefore a factorization of $u$.
By Theo\-rem~\ref{2main}~(c), it thus equals the canonical identification
$\Gr_0 M_- = \Gr_0 M_+$. Hence $\Gr_0 M_- = \Gr_0 M_+$ is a retract
of $N$. Consider a cone of $\Gr_0 M_- \to N$ in $\CC$:
\[ 
\Gr_0 M_- \longto N 
\longto P \longto \Gr_0 M_-[1]
\]
The exact triangle is split in the sense that $P \to \Gr_0 M_-[1]$
is zero. Hence the morphism $N \to P$ admits a right inverse $P \to N$,
unique up to morphisms $P \to \Gr_0 M_-$. There is therefore a unique
right inverse $i: P \to N$ such that its composition with the projection
$p: N \to \Gr_0 M_+ = \Gr_0 M_-$ is zero. The image of $i$ is then 
a kernel of $p$, whose existence is thus established.
This is the canonical complement of $\Gr_0 M_- = \Gr_0 M_+$ in $N$.   

As a retract of $N$, 
the object $P$ belongs to $\CC_{w = 0}$ (condition~\ref{1A}~(1)).
\end{Proof}

\begin{Rem}
Our proof of Corollary~\ref{2E} 
uses the triangulated structure of the category $\CC$.
This can be avoided when the heart 
$\CC_{w = 0}$ is pseudo-Abelian. Namely, consider the composition
\[
p: N \longto Gr_0 M_+ = \Gr_0 M_- \longto N \; . 
\]
It is an idempotent whose image 
is identified with $\Gr_0 M_- = \Gr_0 M_+$. 
Since $\CC_{w = 0}$ is pseudo-Abelian, the morphism $p$ 
also admits a kernel. 
\end{Rem}

For future use, we check the compatibility of Assumption~\ref{2D}
with tensor products. 
Assume therefore that our category $\CC$ is tensor triangulated.
Thus, a bilinear bifunctor
\[
\otimes : \CC \times \CC \longto \CC
\]
is given, and it is assumed to be triangulated in both arguments. 
Assume also that the weight structure $w$ is compatible with $\otimes$,
i.e., that 
\[
\CC_{w \le 0} \otimes \CC_{w \le 0} \subset \CC_{w \le 0} \quad \text{and}
\quad \CC_{w \ge 0} \otimes \CC_{w \ge 0} \subset \CC_{w \ge 0} \; .
\]
It follows that the heart $\CC_{w = 0}$ is a tensor category.
Now fix a second set of data as above.
\begin{enumerate}
\item[(1')] A morphism $u': M_-' \to M_+'$ in $\CC$ between 
$M_-' \in \CC_{w \le 0}$ and $M_+' \in \CC_{w \ge 0} \ $.
\item[(2')] An exact triangle
\[
C' \stackrel{v_-'}{\longto} M_-' \stackrel{u'}{\longto}
M_+' \stackrel{v_+'}{\longto} C'[1] \; .
\]
\end{enumerate}
Fix an exact triangle
\[
D \longto M_- \otimes M_-' \stackrel{u \otimes u'}{\longto}
M_+ \otimes M_+' \longto D[1] \; .
\]

\begin{Prop} \label{2F}
If $C$ and $C'$ are without weights $-1$ and $0$, then so is $D$.
In other words, the validity of Assumption~\ref{2D} for $u$ and $u'$
implies the validity of Assumption~\ref{2D} for $u \otimes u'$.
\end{Prop}

\begin{Proof}
We leave it to the reader to first construct an exact triangle
\[
D \longto M_+ \otimes C' \stackrel{v_+ \otimes v_-'}{\longto}
C \otimes M_-'[1] \longto D[1] \; ,
\]
i.e., to show that $D[1]$ is isomorphic to the cone of the morphism
$v_+ \otimes v_-'$. Then, consider the morphisms
\[
\delta_- \otimes v_-' c_-' :
\Gr_0 M_- \otimes C_{\le -2}' \longto C_{\le -2} \otimes M_-'[1]
\]
and
\[
(c_+[1]) v_+ \otimes \delta_+' :
M_+ \otimes C_{\ge 1}' \longto C_{\ge 1} \otimes \Gr_0 M_+'[1]
\] 
(notation as in Theorem~\ref{2main}). They are completed to give
exact triangles
\[
D_{\le -2} \longto \Gr_0 M_- \otimes C_{\le -2}' 
\longto
C_{\le -2} \otimes M_-'[1] \longto D_{\le -2}[1] 
\]
and
\[
D_{\ge 1} \longto M_+ \otimes C_{\ge 1}'
\longto
C_{\ge 1} \otimes \Gr_0 M_+'[1] \longto D_{\le 1}[1] \; .
\]
By compatibility of $w$ and $\otimes$, and by Proposition~\ref{1C}, 
the object $D_{\le -2}$ is of weights $\le -2$,
and $D_{\ge 1}$ is of weights $\ge 1$.
Finally, it remains to construct an exact triangle
\[
D_{\le -2} \longto D \longto D_{\ge 1} \longto D_{\le -2}[1] \; .
\]
We leave this to the reader (hint: use Theorem~\ref{2main}~(c)).
\end{Proof}

\begin{Cor}
Under the hypotheses of Proposition~\ref{2F}, 
the canoni\-cal morphisms
\[
\Gr_0 \bigl( M_- \otimes M_-' \bigr ) \longto \Gr_0 M_- \otimes \Gr_0 M_-'
\]
and 
\[
\Gr_0 M_+ \otimes \Gr_0 M_+' \longto \Gr_0 \bigl( M_+ \otimes M_+' \bigr ) 
\]
are isomorphisms.
\end{Cor}

%%% Local Variables:

%%% mode: latex
%%% TeX-master: "head"
%%% End:

\bigskip
%\include{Sec4}

%%%%%%%%%%%%%%%%%%%%%%%%%%%%%%%%%%%%%%%%%%%%%%%%%%%%%%%%%%%%%%%%%%%%%%%
%
%  Section 4
%
%%%%%%%%%%%%%%%%%%%%%%%%%%%%%%%%%%%%%%%%%%%%%%%%%%%%%%%%%%%%%%%%%%%%%%%

\section{Example: motives for modular forms}
\label{4}

%%%%%%%%%%%%%%%%%%%%%%%%%%%%%%%%

In his article \cite{Scho}, Scholl constructs the Grothendieck motive $M(f)$
for elliptic normalized newforms $f$ of 
fixed level $n$ and weight $w = r + 2$,
for positive integers $n \ge 3$ and $r \ge 1$.
It is a direct factor of a Grothendieck motive, which 
underlies a Chow motive denoted ${}^r_n \CW$ in \loccit \
(this Chow motive depends only on $n$ and $r$). \\

In order to establish the relation of ${}^r_n \CW$ to the theory
of weights,
let us begin by setting up the notation. It is identical to the one 
introduced in \cite{Scho}, up to one exception: 
the letter $M$ used in \loccit \ to denote certain 
sub-schemes of the modular curve will change to
$S$ in order to avoid
confusion with the motivic notation used earlier in the present paper.
Thus, for our fixed $n \ge 3$ and $r \ge 1$, let $S_n \in Sm/\BQ$
denote the modular curve parametrizing elliptic curves
with level $n$ structure, $j: S_n \into \bS_n$ its smooth
compactification, and $S_n^\infty$ the complement of $S_n$ in $\bS_n$. 
Thus, $S_n^\infty$ is of dimension zero.
Write $X_n \to S_n$ for the universal elliptic curve, and 
$\bX_n \to \bS_n$ for the universal generalized elliptic
curve. Thus, $\bX_n$ is smooth and proper over $\BQ$. The $r$-fold
fibre product 
$\bX_n^r := \bX_n \times_{\bS_n} \times \ldots \times_{\bS_n} \bX_n$ 
of $\bX_n$ over $\bS_n$ is singular
for $r \ge 2$,
and can be desingularized canonically \cite[Lemmas~5.4, 5.5]{D}
(see also \cite[Sect.~3]{Scho}). Denote by $\bbX_n^r$ this desingularization
($\bbX_n^r = \bX_n^r$ for $r = 1$).
Write $X_n^r$ for the $r$-fold
fibre product $X_n$ over $S_n$.  
The symmetric group $\FS_r$ acts on $\bX_n^r$ by permutations,
the $r$-th power of the group $\BZ / n \BZ$ by translations, and the
$r$-th power of the
group $\mu_2$ by inversion in the fibres. Altogether \cite[Sect.~1.1.1]{Scho}, 
this gives a canonical action of the semi-direct product
\[
\Gamma_r := \bigl( (\BZ / n \BZ)^2 \rtimes \mu_2 \bigr)^r \rtimes \FS_r
\]
by automorphisms on $\bX_n^r$. By the canonical nature of the desingularisation,
this extends to an action of $\Gamma_r$ by automorphisms on $\bbX_n^r$.
Of course, this action respects the open sub-scheme $X_n^r$ of $\bbX_n^r$. \\

As in \cite[Sect.~1.1.2]{Scho}, let $\varepsilon: \Gamma_r \to \{ \pm 1 \}$
be the morphism which is trivial on $(\BZ / n \BZ)^{2r}$, is the
product map on $\mu_2^r$, and is the sign character on $\FS_r$.

\begin{Def} \label{4A}
(a)~Let $F$ denote the $\BZ$-algebra $\BZ[1 / (2n \cdot r!)]$. \\[0.1cm]
(b)~Let $e$ denote the idempotent in the group ring
$F[\Gamma_r]$ associated to $\varepsilon \ $:
\[
e := \frac{1}{(2n^2)^r \cdot r!} 
\sum_{\gamma \in \Gamma_r} \varepsilon(\gamma)^{-1} \cdot \gamma 
= \frac{1}{(2n^2)^r \cdot r!} 
\sum_{\gamma \in \Gamma_r} \varepsilon(\gamma) \cdot \gamma  
\]
(observe that $\varepsilon^{-1} = \varepsilon$).
\end{Def}

Let $M$ be an object of an $F$-linear pseudo-Abelian category.
If $M$ comes equipped with an action of the group $\Gamma_r$,
let us agree to denote by $M^e$ the direct factor of $M$ on which
$\Gamma_r$ acts via $\varepsilon$ (in other words, the image of $e$). 
Let us define the object ${}^r_n \CW$ as in \cite[1.2.2]{Scho}.  

\begin{Def} \label{4B} 
Denote by ${}^r_n \CW := \Mgm \bigl( \bbX_n^r \bigr)^e \in \DeffgQM_F$
the image of the idempotent $e$ on $\Mgm \bigl( \bbX_n^r \bigr)$.
\end{Def}

Given that $\bbX_n^r$ is smooth and proper, we see that ${}^r_n \CW$
is an effective Chow motive over $\BQ$.
As above, denote by $\Mgm \bigl( X_n^r \bigr)^e$ and
$\Mcgm \bigl( X_n^r \bigr)^e$ the images of $e$ on 
$\Mgm \bigl( X_n^r \bigr)$ and $\Mcgm \bigl( X_n^r \bigr)$, respectively.
The following can be seen as a translation into the language of 
geometrical motives
of the detailed analysis from
\cite[Sect.~2, 3]{Scho} of the geometry of the boundary of $\bbX_n^r$.

\forget{
We want to apply the results of Section~\ref{3} to $X = X_n^r$,
and to the idempotent $e$.
Thus, we consider the following data.
\begin{enumerate}
\item[(1)] The morphism 
$u: \Mgm \bigl( X_n^r \bigr)^e \to \Mcgm \bigl( X_n^r \bigr)^e$ in 
$\DeffgQM_F \ $. 
\item[(2)] The exact triangle
\[
\dMgm \bigl( X_n^r \bigr)^e \stackrel{v_-}{\longto} 
\Mgm \bigl( X_n^r \bigr)^e 
\stackrel{u}{\longto} \Mcgm \bigl( X_n^r \bigr)^e 
\longto \dMgm \bigl(X_n^r \bigr)^e[1] 
\]
in $\DeffgQM_F \ $.
\end{enumerate}

Application of Corollary~\ref{3B} then implies the following result.

\begin{Cor} \label{4d}
The morphism 
\[
j_n^r : 
\Mgm \bigl( X_n^r \bigr)^e \longto \Mgm \bigl( \bbX_n^r \bigr)^e = {}^r_n \CW
\]
induced by the open immersion $j_n^r$ of $X_n^r$ into $\bbX_n^r$
factors canonically through a monomorphism
\[
\Gr_0 j_n^r : \Gr_0 \Mgm \bigl( X_n^r \bigr)^e \longinto {}^r_n \CW \; .
\]
of effective Chow motives.
This monomorphism is canonically split. 
\end{Cor}

\begin{Prop} \label{4e}
The split monomorphism $\Gr_0 j_n^r$ induces an isomorphism of Grothendieck
motives underlying the Chow motives
$\Gr_0 \Mgm \bigl( X_n^r \bigr)^e$ and ${}^r_n \CW$.
\end{Prop} 

\begin{Proof}
\end{Proof}

\bigskip

\begin{Proofof}{Theorem~\ref{4c}}

\end{Proofof}

\begin{Rem} 
Using the detailed analysis from
\cite[Sect.~2, 3]{Scho} of the geometry of 
the desingularizations $\bbX_n^r$,
one can show that 
\[
\Gr_0 j_n^r : \Gr_0 \Mgm \bigl( X_n^r \bigr)^e \longto {}^r_n \CW 
\]
is actually an isomorphism of Chow motives, for any $r \ge 1$.
}

\begin{Thm} \label{4Ca}
(a)~The motive $\Mcgm \bigl( X_n^r \bigr)^e$ is without weight $1$. 
In particular, the object $\Gr_0 \Mcgm \bigl( X_n^r \bigr)^e$ is defined. 
\\[0.1cm]
(b)~The restriction 
\[
j_n^{r,*} : 
{}^r_n \CW = \Mgm \bigl( \bbX_n^r \bigr)^e \longto
\Mcgm \bigl( X_n^r \bigr)^e  
\]
induced by the open immersion $j_n^r$ of $X_n^r$ into $\bbX_n^r$
factors canonically through an isomorphism
\[
\Gr_0 j_n^{r,*} : {}^r_n \CW \isoto \Gr_0 \Mcgm \bigl( X_n^r \bigr)^e \; .
\]
(c)~There is an exact triangle in $\DeffgQM_F$
\[
C_r \stackrel{\delta_+}{\longto} 
\Gr_0 \Mcgm(X_n^r)^e \stackrel{i_0}{\longto} 
\Mcgm(X_n^r)^e \stackrel{p_+}{\longto} 
C_r[1] \; ,
\]  
where
\[
C_r = \Mgm \bigl( S_n^\infty \bigr) [r]
\]
is pure of weight $r$. 
The exact triangle is canonical up to 
a replacement of the triple of morphisms $(\delta_+,i_0,p_+)$ by
$((-1)^r \delta_+,i_0,(-1)^r p_+)$.
\end{Thm}

Before giving the proof of Theorem~\ref{4Ca},
let us list some of its consequences.
First, duality for smooth schemes \cite[Thm.~4.3.7~3]{V}
implies the following. 

\begin{Cor} \label{4C}
(a)~The motive $\Mgm \bigl( X_n^r \bigr)^e$ is without weight $-1$. 
In particular, the object $\Gr_0 \Mgm \bigl( X_n^r \bigr)^e$ is defined. 
\\[0.1cm]
(b)~The morphism 
\[
j_n^r : 
\Mgm \bigl( X_n^r \bigr)^e \longto \Mgm \bigl( \bbX_n^r \bigr)^e = {}^r_n \CW
\]
factors canonically through an isomorphism
\[
\Gr_0 j_n^r : \Gr_0 \Mgm \bigl( X_n^r \bigr)^e \isoto {}^r_n \CW \; .
\]
(c)~There is an exact triangle in $\DeffgQM_F$
\[
C_{-(r+1)} \stackrel{\iota_-}{\longto} \Mgm(X_n^r)^e 
\stackrel{\pi_0}{\longto} \Gr_0 \Mgm(X_n^r)^e 
\stackrel{\delta_-}{\longto} C_{-(r+1)}[1] \; ,
\]  
where
\[
C_{-(r+1)} = \Mgm \bigl( S_n^\infty \bigr) (r+1)[r+1]
\]
is pure of weight $-(r+1)$. 
The exact triangle is canonical up to 
a replacement of the triple of morphisms $(\iota_-,\pi_0,\delta_-)$ by
$((-1)^r \iota_-,\pi_0,(-1)^r \delta_-)$.
\end{Cor}

\begin{Rem} \label{4D}
(a)~It is well known that the motive $\Mgm \bigl( S_n^\infty \bigr)$
is isomorphic to a finite sum of copies of $\Mgm (\Spec \BQ(\mu_n) )$. \\[0.1cm]
(b)~As the proof will show, 
Theorem~\ref{4Ca}~(c) and Corollary~\ref{4C}~(c) remain true for $r =0$
if one replaces $\Gr_0 \Mcgm(X_n^r)^e$ and $\Gr_0 \Mgm(X_n^r)^e$
by $\Mgm (\bS_n)$.
\end{Rem}

Next, note that Theorem~\ref{4Ca}~(b) and Corollary~\ref{4C}~(b) together 
imply the following.

\begin{Cor} \label{4Da}
The canonical morphism $\Mgm \bigl( X_n^r \bigr)^e \to 
\Mcgm \bigl( X_n^r \bigr)^e$ factors canonically through an isomorphism
\[ 
\Gr_0 \Mgm \bigl( X_n^r \bigr)^e \isoto
\Gr_0 \Mcgm \bigl( X_n^r \bigr)^e \; .
\]
\end{Cor}

For any object $M$ of $\DeffgQM_F \ $, define motivic cohomology 
\[
H_{\CM}^p \bigl( M,F(q) \bigr) :=  
\Hom_{\DeffgQM_F} \bigl( M,\BZ(q)[p] \bigr) \; .
\]
When $M = \Mgm(Y)$ for a scheme $Y \in Sm / \BQ$, 
this gives motivic cohomology
$H_{\CM}^p \bigl( Y,\BZ(q) \bigr)$ of $Y$,
tensored with $F = \BZ[1 / (2n \cdot r!)]$. Thus for example,
\[
H_{\CM}^{r+1} \bigl( C_{-(r+1)},F(r+1) \bigr) =
H_{\CM}^0 \bigl( S_n^\infty,\BZ(0) \bigr) \otimes_\BZ F \; .
\]
Similarly,
\[
H_{\CM}^{r+2} \bigl( C_{-(r+1)},F(r+\ell+2) \bigr) =
H_{\CM}^1 \bigl( S_n^\infty,\BZ(\ell+1) \bigr) 
\otimes_\BZ F 
\]
for any integer $\ell$.
We get the following refinement of \cite[Cor.~1.4.1]{Scho}.

\begin{Cor} \label{4F}
Let $\ell \ge 0$ be a second integer. Then 
the kernel of the morphism 
\[
H_{\CM}(\iota_-):
\bigl( H_{\CM}^{r+2} \bigl( X_n^r,\BZ(r+\ell+2) \bigr) \otimes_\BZ F \bigr)^e
\longto H_{\CM}^1 \bigl( S_n^\infty,\BZ(\ell+1) \bigr) 
\otimes_\BZ F 
\]
equals  
$H_{\CM}^{r+2} \bigl( \Gr_0 \Mgm(X_n^r)^e,F(r+\ell+2) \bigr)$.
\end{Cor}

\begin{Proof}
This follows from the exact triangle of Corollary~\ref{4C}~(c), and from the
vanishing of $H_{\CM}^0 \bigl( S_n^\infty,\BZ(\ell+1) \bigr)$
(since $\ell +1 \ge 1$).
\end{Proof}

Recall that following ideas of Beilinson \cite{B},
this result can be employed as follows: 
using the \emph{Eisenstein symbol}
defined in \loccit \, one constructs elements in
$\bigl( H_{\CM}^{r+2} \bigl( X_n^r,\BZ(r+\ell+2) \bigr) 
\otimes_\BZ \BQ \bigr)^e$.
By Corollary~\ref{4F},
linear combinations of such elements vanishing under $H_{\CM}(\iota_-)$
lie in the sub-$\BQ$-vector space 
$H_{\CM}^{r+2} \bigl( \Gr_0 \Mgm(X_n^r)^e,\BQ(r+\ell+2) \bigr)$.
It can then be shown that there are sufficiently many such linear
combinations, in the sense that their images under the regulator 
generate \emph{Deligne cohomology} \cite[Sect.~2 and 4]{Schn}
$H_{\CD}^{r+2} \bigl( \Gr_0 \Mgm(X_n^r)^e_{/ \BR},\BR(r+\ell+2) \bigr)$.
Furthermore, the $\BQ$-span of these images
has the relation to the leading coefficient 
at $s = - \ell$ of the $L$-function of 
$\Gr_0 \Mgm(X_n^r)^e$, predicted by Beilinson's conjecture. 
For details, see forthcoming work of Scholl \cite{Scho2}. \\

\begin{Proofof}{Theorem~\ref{4Ca}}
Denote by $\bbX_n^{r,\infty}$ the complement of the smooth scheme
$X_n^r$ in the smooth and proper scheme $\bbX_n^r$.
Localization for the motive with compact support \cite[Prop.~4.1.5]{V} 
shows that there is a canonical exact triangle
\[
\Mgm \bigl( \bbX_n^r \bigr)^e \stackrel{j_n^{r,*}}{\longto} 
\Mcgm(X_n^r)^e \longto
\Mgm \bigl( \bbX_n^{r,\infty} \bigr)^e [1] \longto
\Mgm \bigl( \bbX_n^r \bigr)^e [1] \; .
\]
Following the strategy from \cite{Scho},
we shall show the following claim.
\[
(C) \quad \quad \quad \quad 
\Mgm \bigl( \bbX_n^{r,\infty} \bigr)^e = \Mcgm \bigl( \bbX_n^{r,\infty} \bigr)^e
\isoto \Mcgm \bigl( S_n^\infty \bigr)[r] = \Mgm \bigl( S_n^\infty \bigr)[r] 
\]
canonically up to a sign $(-1)^r$.
In particular, the motive $\Mgm \bigl( \bbX_n^{r,\infty} \bigr)^e$
is pure of weight $r \ge 1$.
Claim~(C) implies that the above exact triangle 
is a weight filtration of $\Mcgm(X_n^r)^e$
avoiding weight $1$. Furthermore, the restriction $j_n^{r,*}$ identifies
$\Mgm \bigl( \bbX_n^r \bigr)^e$ with $\Gr_0 \Mcgm(X_n^r)^e$.

To show claim $(C)$, observe first that 
the motivic version of \cite[Statement~1.3.0]{Scho}
remains valid: for any $S \in Sm / \BQ$,
there is a decomposition in $\DeffgQM_F$ (in fact, already
in $\DeffgQM_{\BZ[1/2]}$)
\[
\Mcgm (\BG_m \times_\BQ S) \cong \Mcgm (S)(1)[2] \oplus \Mcgm (S)[1] \; ,
\]
such that inversion $x \mapsto x^{-1}$ on $\BG_m$ acts on the first factor
by $+1$, and on the second by $-1$. The projection onto the first factor 
is canonical, and the projection onto the second factor is canonical
up to a sign. Proof: localization \cite[Prop.~4.1.5]{V}
for the inclusion of $\BG_m$ into the projective line; the choice
of the second projection is equivalent to the choice of one of the
residue morphisms to $0$ or $\infty$.
Fix one of the two choices of projection
\[
\pi_- : \Mcgm (\BG_m) \longonto \Mcgm (\Spec \BQ)[1] 
\]
(and use the same notation for 
$\Mcgm (\BG_m \times_\BQ S) \onto \Mcgm (S)[1]$ obtained by base change
via $S$). Then to define the morphism
\[
j^{0,*} : \Mcgm \bigl( \bbX_n^{r,\infty} \bigr)^e \longto
\Mcgm \bigl( S_n^\infty \bigr)[r] \; , 
\]
consider the following.
\begin{enumerate}
\item[(i)] The 
intersection $\bbX_n^{r,\infty,\reg}$
of $\bbX_n^{r,\infty}$ 
with the non-singular part $\bX_n^{r,\reg}$ of $\bX_n^r$. Explanation: by
\cite[Thm.~3.1.0~ii)]{Scho}, the desingularization $\bbX_n^r \onto \bX_n^r$
is an isomorphism over $\bX_n^{r,\reg}$. Thus, $\bbX_n^{r,\infty,\reg}$
is an open sub-scheme of $\bbX_n^{r,\infty}$,
\item[(ii)] the 
neutral component $\bbX_n^{r,\infty,0}$
of $\bbX_n^{r,\infty,\reg}$, i.e., its 
intersection with the N\'eron model of $X_n^r$.
Thus, $\bbX_n^{r,\infty,0}$ is an open sub-scheme of $\bbX_n^{r,\infty,\reg}$,
\item[(iii)] the identification of $\bbX_n^{r,\infty,0}$
with the $r$-th power (over the base 
$S_n^\infty$) of the neutral
component $\bbX_n^{1,\infty,0}$ for $r = 1$.
The latter can be identified with $\BG_m \times_\BQ S_n^\infty$,
canonically up to an automorphism $x \mapsto x^{- 1}$. Hence 
$\bbX_n^{r,\infty,0} \cong \BG_m^r \times_\BQ S_n^\infty$,
canonically up to an automorphism 
\[
(x_1,\ldots,x_r) \longmapsto (x_1,\ldots,x_r)^{- 1} \; .
\]
\end{enumerate}
The three steps (i)--(iii) give an open immersion 
\[
j^0 : \BG_m^r \times_\BQ S_n^\infty \longinto \bbX_n^{r,\infty} \; ,
\]
which by contravariance of $\Mcgm$ induces a morphism
\[
j^{0,*} : \Mcgm \bigl( \bbX_n^{r,\infty} \bigr)^e \longto
\Mcgm \bigl( \BG_m^r \times_\BQ S_n^\infty \bigr) \; .
\]
Its composition with the $r$-th power of $\pi_-$ gives the desired
morphism
\[
\Mcgm \bigl( \bbX_n^{r,\infty} \bigr)^e \longto
\Mcgm \bigl( S_n^\infty \bigr) [r] \; ,
\] 
equally denoted $j^{0,*}$, and canonical up to a sign $(-1)^r$.
In the context of twisted Poincar\'e duality theories $(H^*,H_*)$,
Scholl's main technical result 
\cite[Thm.~1.3.3]{Scho} is equivalent to stating that the morphism
induced by $j^{0,*}$ on the level of the theory $H_*$ is an isomorphism.

Our observation is simply that the same proof as the one given 
in \loccit \ runs through, with $(H^*,H_*)$ replaced by 
$(\Mgm,\Mcgm)$. 

More precisely, \cite[Lemma~1.3.1]{Scho} holds for $\Mcgm$, and hence
the proof of \cite[Prop.~2.4.1]{Scho} runs through for $\Mcgm$. The
latter result implies that the schemes occurring in a suitable stratification
of the complement of $\bbX_n^{r,\infty,\reg}$ in $\bbX_n^{r,\infty}$
all have trivial $M_{gm}^{c,e}$ 
(cmp.\ \cite[proof of Thm.~3.1.0~ii)]{Scho}).
This shows that step (i) induces an isomorphism
\[
\Mcgm \bigl( \bbX_n^{r,\infty} \bigr)^e \isoto
\Mcgm \bigl( \bbX_n^{r,\infty,\reg} \bigr)^e \; .
\]
To deal with step (ii), one shows that the group $\Gamma_r$ acts
transitively on the set of components of the singular part of
$\bbX_n^{r,\infty,\reg}$, and that the stabilizer
of each component admits a 
subgroup of order two acting trivially on the component,
but having trivial intersection with the kernel of $\varepsilon$ 
(cmp.\ \cite[proof of Thm.~3.1.0~iii)]{Scho}). This shows first that
the singular part of $\bbX_n^{r,\infty,\reg}$ does not contribute to
$M_{gm}^{c,e}$, and then that 
\[
\Mcgm \bigl( \bbX_n^{r,\infty,\reg} \bigr)^e  \longto
\Mcgm \bigl( \bbX_n^{r,\infty,0} \bigr)^{e'}
\]
is an isomorphism, where $e'$ denotes the projection onto the 
eigenspace for the restriction of $\varepsilon$ to the subgroup
$\mu_2^r \rtimes \FS_r$ of $\Gamma_r$. To conclude, we 
apply the motivic version of
\cite[Lemma~1.3.1]{Scho} to see that $\pi_-^{\otimes r}$ induces an 
isomorphism
$\Mcgm \bigl( \BG_m^r \times_\BQ S_n^\infty \bigr)^{e'} \isoto 
\Mcgm \bigl( S_n^\infty \bigr)[r]$. 
\end{Proofof}

\begin{Rem} \label{4H}
(a)~The proof of Theorem~\ref{4Ca} shows that the open immersion 
of the non-singular part $\bX_n^{r,\reg}$ of $\bX_n^r$ into $\bbX_n^r$
induces an isomorphism
\[
\Mgm \bigl( \bbX_n^r \bigr)^e \isoto
\Mcgm \bigl( \bX_n^{r,\reg} \bigr)^e \; .
\]
In particular, $\Mcgm \bigl( \bX_n^{r,\reg} \bigr)^e$ is a Chow motive.
By duality for smooth schemes \cite[Thm.~4.3.7~3]{V},
\[
\Mgm \bigl( \bX_n^{r,\reg} \bigr)^e \longto
\Mgm \bigl( \bbX_n^r \bigr)^e 
\]
is an isomorphism, too, and hence so is the canonical morphism
\[
\Mgm \bigl( \bX_n^{r,\reg} \bigr)^e \longto
\Mcgm \bigl( \bX_n^{r,\reg} \bigr)^e \; .
\]
The construction of the motives 
$M(f)$ can therefore also be done using the smooth non-proper
scheme $\bX_n^{r,\reg}$
instead of $\bbX_n^r$. \\[0.1cm]
(b)~A slightly closer look at the proof of \cite[Thm.~3.1.0~ii)]{Scho} 
reveals that the open immersion 
of $\bX_n^{r,\reg}$ into $\bX_n^r$
also induces an isomorphism
\[
\Mgm \bigl( \bX_n^r \bigr)^e \isoto
\Mcgm \bigl( \bX_n^{r,\reg} \bigr)^e \; .
\]
By (a), $\Mgm \bigl( \bX_n^r \bigr)^e = \Mcgm \bigl( \bX_n^r \bigr)^e$ 
is a Chow motive.
The construction of the motives 
$M(f)$ can therefore also be done using the non-smooth (for $r \ge 2$)
proper scheme $\bX_n^r$ instead of $\bbX_n^r$. 
\end{Rem}

%%%%%%%%%%%%%%%%%%%%%%%%%%%%%%%%

%%% Local Variables:
%%% mode: latex
%%% TeX-master: "head"
%%% End:

\bigskip
%\include{Sec3}

%%%%%%%%%%%%%%%%%%%%%%%%%%%%%%%%%%%%%%%%%%%%%%%%%%%%%%%%%%%%%%%%%%%%%%%
%
%  Section 3
%
%%%%%%%%%%%%%%%%%%%%%%%%%%%%%%%%%%%%%%%%%%%%%%%%%%%%%%%%%%%%%%%%%%%%%%%

\section{Weights, boundary motive and interior motive}
\label{3}

%%%%%%%%%%%%%%%%%%%%%%%%%%%%%%%%

%%%%%%%%%%%%%%%%%%%%%%%%%%%%%%%%

This section contains our main result (Theorem~\ref{Main}). 
We list its main consequences, and define in particular
the motivic analogue of (certain direct factors of) interior
cohomology (Definition~\ref{3def}). 
Throughout, we assume $k$ to admit resolution of singularities. \\ 

Let us fix $X \in Sm/k$.
The \emph{boundary motive} $\dMgm(X)$ of $X$ \cite[Def.~2.1]{W1}
fits into a canonical exact triangle
\[
(\ast) \quad\quad
\dMgm(X) \longto \Mgm(X) \longto \Mcgm(X) \longto \dMgm(X)[1]
\]
in $\DeffgM$. The algebra of \emph{finite correspondences} 
$c(X,X)$ acts on $\Mgm(X)$ \cite[p.~190]{V}. Denote by 
${}^t c(X,X)$ the transposed algebra: a cycle $\FZ$ on $X \times_k X$
lies in ${}^t c(X,X)$ if and only if ${}^t \FZ \in c(X,X)$.
The definition of composition of correspondences \loccit \
shows that the intersection $c(X,X) \cap {}^t c(X,X)$ acts on $\Mcgm(X)$. 

\begin{Def}  \label{3a}
(a)~Define the algebra $c_{1,2}(X,X)$ as the intersection of the algebras
$c(X,X)$ and ${}^t c(X,X)$. As an Abelian group, $c_{1,2}(X,X)$ is thus
free on the symbols $(Z)$, where $Z$ runs through the integral closed 
sub-schemes of $X \times_k X$, such that both projections to the
components $X$ are finite on $Z$, and map $Z$ surjectively to a connected
component of $X$. Multiplication in $c_{1,2}(X,X)$ is defined by composition
of correspondences as in \cite[p.~190]{V}. \\[0.1cm]
(b)~Denote by ${}^t$ the canonical anti-involution on $c_{1,2}(X,X)$
mapping a cycle $\FZ$ to ${}^t \FZ$.
\end{Def}

It results directly from the definitions that the algebra $c_{1,2}(X,X)$
acts on the triangle $(\ast)$
in the sense that it acts on the three objects, and the morphisms 
are $c_{1,2}(X,X)$-equivariant. Denote by $\bar{c}_{1,2}(X,X)$ 
the quotient of $c_{1,2}(X,X)$ by the kernel of this action. 
Fix a commutative flat $\BZ$-algebra $F$, and an idempotent $e$
in $\bar{c}_{1,2}(X,X) \otimes_\BZ F$. 
Denote by $\Mgm(X)^e$, $\Mcgm(X)^e$ and $\dMgm(X)^e$ the images of $e$
on $\Mgm(X)$, $\Mcgm(X)$ and $\dMgm(X)$, respectively, considered as
objects of the category $\DeffgM_F$. 
We are ready to set up the data (1), (2) considered in Section~\ref{2},
for $\CC := \DeffgM_F$.
\begin{enumerate}
\item[(1)] The morphism $u$ is the morphism
$\Mgm(X)^e \to \Mcgm(X)^e$. By Corollary~\ref{1E}
and condition \ref{1A}~(1), the object $\Mgm(X)^e$
belongs indeed to $\DeffgM_{F,w \le 0}$, and 
$\Mcgm(X)^e$ to $\DeffgM_{F,w \ge 0} \ $.
\item[(2)] Our choice of cone of $u$ is $\dMgm(X)^e[1]$, together with
the exact triangle
\[
\dMgm(X)^e \stackrel{v_-}{\longto} \Mgm(X)^e \stackrel{u}{\longto}
\Mcgm(X)^e \stackrel{v_+}{\longto} \dMgm(X)^e[1] 
\]
in $\DeffgM_F$ induced by $(\ast)$.
\end{enumerate}

Observe that our data (1), (2) are stable under the natural action
of 
\[
GCen_{\bar{c}_{1,2}(X,X)}(e) := 
\big{ \{ } z \in \bar{c}_{1,2}(X,X) \otimes_\BZ F \; , \;
ze = eze \big{ \} } \; .
\]
In particular, they are stable under the action of the
centralizer $Cen_{\bar{c}_{1,2}(X,X)}(e)$
of $e$ in $\bar{c}_{1,2}(X,X) \otimes_\BZ F$.
In this context, Assumption~\ref{2D} reads as follows.

\begin{Ass} \label{3A}
The direct factor $\dMgm(X)^e$ of the boundary motive of $X$
is without weights $-1$ and $0$.
\end{Ass}

Thus, we may and do fix a weight filtration 
\[ 
C_{\le -2} \stackrel{c_-}{\longto} \dMgm(X)^e 
\stackrel{c_+}{\longto} C_{\ge 1} \stackrel{\delta_C}{\longto} C_{\le -2}[1]
\]
avoiding weights $-1$ and $0$. Theorem~\ref{2main}, Corollary~\ref{2E}
and the adjunction property from Proposition~\ref{2Cc}
then give the following.

\begin{Thm} \label{Main}
Fix the data (1), (2), and suppose Assumption~\ref{3A}. \\[0.1cm]
(a)~The motive $\Mgm(X)^e$ is without weight $-1$,
and the motive $\Mcgm(X)^e$ is without weight $1$. 
In particular, the effective Chow motives $\Gr_0 \Mgm(X)^e$ 
and $\Gr_0 \Mcgm(X)^e$ are defined,
and they carry a natural action of $GCen_{\bar{c}_{1,2}(X,X)}(e)$. \\[0.1cm]
(b)~There are canonical exact triangles 
\[
(3) \quad\quad
C_{\le -2} \stackrel{v_- c_-}{\longto} \Mgm(X)^e 
\stackrel{\pi_0}{\longto} \Gr_0 \Mgm(X)^e 
\stackrel{\delta_-}{\longto} C_{\le -2}[1]
\]
and
\[
(4) \quad\quad
C_{\ge 1} \stackrel{\delta_+}{\longto} \Gr_0 \Mcgm(X)^e 
\stackrel{i_0}{\longto} \Mcgm(X)^e 
\stackrel{(c_+[1]) v_+}{\longto} C_{\ge 1}[1] \; ,
\]
which are stable under the natural action
of $GCen_{\bar{c}_{1,2}(X,X)}(e)$. \\[0.1cm]
(c)~There is a canonical isomorphism 
$\Gr_0 \Mgm(X)^e \isoto \Gr_0 \Mcgm(X)^e$ in $\CHeffM_F$. 
As a morphism, it is uniquely determined by the property
of making the diagram
\[
\vcenter{\xymatrix@R-10pt{
        \Mgm(X)^e \ar[r]^-{u} \ar[d]_{\pi_0} &
        \Mcgm(X)^e \\
        \Gr_0 \Mgm(X)^e \ar[r] &
        \Gr_0 \Mcgm(X)^e \ar[u]_{i_0}
\\}}
\]
commute; in particular, it is $GCen_{\bar{c}_{1,2}(X,X)}(e)$-equivariant. 
Its inverse makes the diagram
\[
\vcenter{\xymatrix@R-10pt{
        C_{\ge 1} \ar[r]^-{\delta_C} \ar[d]_{\delta_+} &
        C_{\le -2}[1] \\
        \Gr_0 \Mcgm(X)^e \ar[r] &
        \Gr_0 \Mgm(X)^e \ar[u]_{\delta_-}
\\}}
\]
commute. \\[0.1cm]
(d)~Let $N \in \CHM_F$ be a Chow motive. Then $\pi_0$ and $i_0$ induce
isomorphisms
\[
\Hom_{\CHM_F} \bigl( \Gr_0 \Mgm(X)^e , N \bigr) \isoto
\Hom_{\DgM_F} \bigl( \Mgm(X)^e , N \bigr)
\]
and
\[
\Hom_{\CHM_F} \bigl( N , \Gr_0 \Mcgm(X)^e \bigr) \isoto
\Hom_{\DgM_F} \bigl( N , \Mcgm(X)^e \bigr) \; .
\]
(e)~Let $\Mgm(X)^e \to N \to \Mcgm(X)^e$ be a factorization of $u$ 
through a Chow motive
$N \in \CHM_F$. 
Then $\Gr_0 \Mgm(X)^e = \Gr_0 \Mcgm(X)^e$ is canonically
a direct factor of $N$, with a canonical direct complement.
\end{Thm} 

We explicitly mention the following immediate consequence of 
Theorem~\ref{Main}. 

\begin{Cor} \label{3Aa}
Fix $X$ and $e$, and suppose that $\dMgm(X)^e = 0$, i.e., that
\[
u: \Mgm(X)^e \isoto \Mcgm(X)^e \; .
\]
Then $\Mgm(X)^e \cong \Mcgm(X)^e$ are effective Chow motives.
\end{Cor}

Of course, this also follows from Corollary~\ref{1F}. 
The author knows of no proof of Corollary~\ref{3Aa} ``avoiding weights''
when $e \ne 1$. (For $e = 1$, we leave it to the reader to show 
(using for example \cite[Cor.~4.2.5]{V}) that 
the assumption $\dMgm(X) = 0$ is equivalent to $X$ being proper.) 

\begin{Rem} \label{3Ab}
It is not difficult to see that Assumption~\ref{3A} is actually
implied by parts (a) and (c) of Theorem~\ref{Main}.
\end{Rem}

Henceforth, we identify $\Gr_0 \Mgm(X)^e$ and $\Gr_0 \Mcgm(X)^e$
via the canonical isomorphism of Theorem~\ref{Main}~(c).

\begin{Cor} \label{3B} 
In the situation considered in Theorem~\ref{Main},
let $\Xp$ be any smooth compactification of $X$. Then
$\Gr_0 \Mgm(X)^e$ is canoni\-cally
a direct factor of the Chow motive
$\Mgm(\Xp)$, with a canonical direct complement.
\end{Cor}

\begin{Proof}
Indeed, the morphism $u$ factors canonically through $\Mgm(\Xp)$:
\[
\Mgm(X)^e \longinto \Mgm(X) \longto \Mgm(\Xp)
\longto \Mcgm(X) \longonto \Mcgm(X)^e \; .
\]
Hence we may apply Theorem~\ref{Main}~(e).
\end{Proof}

Recall that the category $\CHM_F$ is
pseudo-Abelian. Thus, the construction of a sub-motive of $\Mgm(\Xp)$ 
does not \emph{a priori}
necessitate the \emph{identification}, but only the 
\emph{existence} of a complement. In our
situation, Corollary~\ref{3B} states that
the complement of $\Gr_0 \Mgm(X)^e$ \emph{is} canonical.
This shows
that Assumption~\ref{3A} is indeed rather restrictive,
an observation confirmed by part~(c) of the
following results on the Hodge theoretic and $\ell$-adic realizations
(\cite[Sect.~2 and Corrigendum]{H}; see \cite[Sect.~1.5]{DG}
for a simplification of this approach).
They can be seen as applications 
of the \emph{cohomological weight spectral sequence}
\cite[Thm.~2.4.1, Rem.~2.4.2]{Bo} in a very special case.

\begin{Thm} \label{3C}
Keep the situation considered in Theorem~\ref{Main}. 
Assume that $k$ can be embedded into the field $\BC$ of complex numbers.
Fix one such embedding.
Let $H^\ast$ be the Hodge theoretic realization
\cite[Cor.~2.3.5 and Corrigendum]{H},
followed by the canonical cohomology functor,
i.e., the functor on $\DeffgM_F$
given by Betti cohomology of the topological space of 
$\BC$-valued points, tensored with $\BQ \otimes_\BZ F$, and
with its natural mixed Hodge structure. 
Let $n \in \BN$. \\[0.1cm]
(a)~The morphisms $\pi_0$ and
$i_0$ induce isomorphisms
\[
H^n \bigl( \Gr_0 \Mgm(X)^e \bigr) \isoto 
W_n H^n \bigl( \Mgm(X)^e \bigr) = 
\bigl( W_n H^n (X (\BC),\BQ) \otimes_\BZ F \bigr)^e 
\]
and
\[
\frac{\bigl( H^n_c (X (\BC),\BQ) \otimes_\BZ F \bigr)^e}{\bigl( 
          W_{n-1} H^n_c (X (\BC),\BQ) \otimes_\BZ F \bigr)^e} = 
\frac{H^n (\Mcgm(X)^e)}{W_{n-1} H^n (\Mcgm(X)^e)} \isoto 
H^n \bigl( \Gr_0 \Mgm(X)^e \bigr) .
\]
Here, $W_r$ denotes the $r$-th filtration step
of the weight filtration of a mixed Hodge structure (thus, the weights of
$H^n (X (\BC),\BQ)$ are $\ge n$,
and those of $H^n_c (X (\BC),\BQ)$ are $\le n$). \\[0.1cm]
(b)~The isomorphisms of (a) identify $H^n (\Gr_0 \Mgm(X)^e)$ with the image
of the natural morphism
\[
\bigl( H^n_c (X (\BC),\BQ) \otimes_\BZ F \bigr)^e \longto 
\bigl( H^n (X (\BC),\BQ) \otimes_\BZ F \bigr)^e \; .
\]
(c)~The image of $( H^n_c (X (\BC),\BQ) \otimes_\BZ F)^e$ 
in $( H^n (X (\BC),\BQ) \otimes_\BZ F)^e$ 
equals the lowest weight filtration step $W_n$ of 
$( H^n (X (\BC),\BQ) \otimes_\BZ F)^e$.
\end{Thm}

The reader should be aware that the algebra $\bar{c}_{1,2}(X,X)$ acts
contravariantly on Betti cohomology $H^n (X (\BC),\BQ)$. The same remark
applies of course to the $\ell$-adic realization, which we consider now. 

\begin{Thm} \label{3D}
Keep the situation considered in Theorem~\ref{Main}, and fix a prime 
$\ell > 0$. 
Assume that $k$ is finitely generated over its prime field, 
and of characteristic zero.
Let $H^\ast$ be the $\ell$-adic realization
\cite[Cor.~2.3.4 and Corrigendum]{H},
followed by the canonical cohomology functor,
i.e., the functor on $\DeffgM_F$
given by $\ell$-adic cohomology of the base change to a fixed algebraic
closure $\bar{k}$ of $k$, tensored with $\BQ_\ell \otimes_\BZ F$, and
with its natural action of the absolute Galois 
group $G_k$ of $k$. Let $n \in \BN$. \\[0.1cm]
(a)~The morphisms $\pi_0$ and
$i_0$ induce isomorphisms
\[
H^n \bigl( \Gr_0 \Mgm(X)^e \bigr) \isoto 
W_n H^n \bigl( \Mgm(X)^e \bigr) = 
\bigl( W_n H^n (X_{\bar{k}},\BQ_\ell) \otimes_\BZ F \bigr)^e 
\]
and
\[
\frac{\bigl( H^n_c (X_{\bar{k}},\BQ_\ell) \otimes_\BZ F \bigr)^e}{\bigl( 
          W_{n-1} H^n_c (X_{\bar{k}},\BQ_\ell) \otimes_\BZ F \bigr)^e} = 
\frac{H^n (\Mcgm(X)^e)}{W_{n-1} H^n (\Mcgm(X)^e)} \isoto 
H^n \bigl( \Gr_0 \Mgm(X)^e \bigr) \; .
\]
Here, $W_r$ denotes the $r$-th filtration step
of the weight filtration of a $G_k$-module,
and $X_{\bar{k}}$ denotes the base change $X \otimes_k \, \bar{k}$
of $X$ to $\bar{k}$ (thus, the weights of
$H^n (X_{\bar{k}},\BQ_\ell)$ are $\ge n$,
and those of $H^n_c (X_{\bar{k}},\BQ_\ell)$ are $\le n$). \\[0.1cm]
(b)~The isomorphisms of (a) identify $H^n (\Gr_0 \Mgm(X)^e)$ with the image
of the natural morphism
\[
\bigl( H^n_c (X_{\bar{k}},\BQ_\ell) \otimes_\BZ F \bigr)^e \longto 
\bigl( H^n (X_{\bar{k}},\BQ_\ell) \otimes_\BZ F \bigr)^e \; .
\]
(c)~The image of $( H^n_c (X_{\bar{k}},\BQ_\ell) \otimes_\BZ F)^e$ 
in $( H^n (X_{\bar{k}},\BQ_\ell) \otimes_\BZ F)^e$ 
equals the lowest weight filtration step $W_n$ of 
$( H^n (X_{\bar{k}},\BQ_\ell) \otimes_\BZ F)^e$.
\end{Thm}

\begin{Def} \label{3def}
Fix the data (1), (2), and suppose Assumption~\ref{3A}.
We call $\Gr_0 \Mgm(X)^e$ the \emph{$e$-part of the
interior motive of $X$}.
\end{Def}

This terminology is motivated by parts~(b) of Theorems~\ref{3C} and \ref{3D},
which show that after passage to rational coefficients, the 
realizations of $\Gr_0 \Mgm(X)^e$ are classes of
complexes, whose cohomology equals the part of the
interior cohomology of $X$ fixed by $e$.

\medskip

\begin{Proofof}{Theorems~\ref{3C} and \ref{3D}}
Consider the exact triangle
\[
(3) \quad\quad
C_{\le -2} \longto \Mgm(X)^e 
\stackrel{\pi_0}{\longto} \Gr_0 \Mgm(X)^e 
\longto C_{\le -2}[1]
\]
from Theorem~\ref{Main}~(b). Recall that 
as suggested by the notation, the motive
$C_{\le -2}$ is of weights $\le -2$.
The cohomological functor $H^*$ transforms
it into a long exact sequence
\[
H^{n-1} \bigl( C_{\le -2} \bigr) \longto H^n \bigl( \Gr_0 \Mgm(X)^e \bigr)
\stackrel{H^n(\pi_0)}{\longto} H^n \bigl( \Mgm(X)^e \bigr) \longto
H^n \bigl( C_{\le -2} \bigr) \; .
\]
The essential information we need to use is that
$H^n$ transforms Chow motives into objects which are pure of weight $n$,
for any $n \in \BN$. Since $C_{\le -2}$ admits a filtration whose
cones are Chow motives sitting in degrees $\ge 2$, the object
$H^n(C_{\le -2})$ admits a filtration whose graded pieces 
are of weights $\ge n +2$. Since our coefficients are $\BQ$-vector spaces,
there are no non-trivial morphisms between objects of disjoint weights. 
The above long exact sequence then shows that 
\[
H^n(\pi_0): H^n \bigl( \Gr_0 \Mgm(X)^e \bigr)
\longto H^n \bigl( \Mgm(X)^e \bigr)
\]
is injective, and its image is identical to the part of weight $n$
of $H^n ( \Mgm(X)^e )$. This shows the part of claim (a)
concerning $\pi_0$.

Now recall that the realization functor is compatible with the tensor
structures \cite[Cor.~2.3.5, Cor.~2.3.4]{H}, and sends the Tate motive
$\BZ(1)$ to the dual of the Tate object \cite[Thm.~2.3.3]{H}. 
It follows that it is compatible
with duality. Since on the one hand, $\Mgm(X)$ and $\Mcgm(X)$
are in duality \cite[Thm.~4.3.7~3]{V}, and on the other hand,
the same is true for Betti, resp.\ $\ell$-adic cohomology
and Betti, resp.\ $\ell$-adic cohomology with compact support,
we see that $H^n$ sends the motive with compact support $\Mcgm(X)$
to cohomology with compact support $H^n_c$ of $X$.

Now repeat the above argument for the exact triangle 
\[
(4) \quad\quad
C_{\ge 1} \longto \Gr_0 \Mcgm(X)^e 
\stackrel{i_0}{\longto} \Mcgm(X)^e 
\longto C_{\ge 1}[1]
\]
from Theorem~\ref{Main}~(b).
This shows the remaining part of claim (a).

Claims (b) and (c) follow, once we observe that the composition of
$H^n(i_0)$ and $H^n(\pi_0)$ equals the canonical morphism from
cohomology with compact support to cohomology without support.  
\end{Proofof}

\begin{Rem} 
(a)~The proof of Theorems~\ref{3C} and \ref{3D}
uses the fact that in the respective target categories
(mixed Hodge structures in Theorem~\ref{3C},
Galois representations in Theorems~\ref{3D}), 
there are no non-trivial morphisms between objects of disjoint weights. 
This is true as long as we work with coefficients which are
$\BQ$-vector spaces. Note that recent work of Lecomte 
\cite{L} establishes the existence of a Betti realization
which does not require the passage to $\BQ$-coefficients.
In particular \cite[Thm.~1.1]{L}, for a smooth quasi-projective
variety $X$, it yields the classical singular cohomology of 
the topological space $X(\BC)$. 
We have no statement (and not even a guess) to offer
on the image of $\Gr_0 \Mgm(X)^e$
under the realization of \loccit . \\[0.1cm]
(b)~The author does not know whether for general $Y \in Sm/k$
it is possible (or even reasonable to expect) to find a
complex computing interior (Betti or $\ell$-adic)
cohomology of $Y$, and through which
the natural morphism $R \Gamma_c (Y) \to R \Gamma (Y)$ factors. 
\end{Rem}

\begin{Rem}
When $k$ is a number field,  
Theorems~\ref{3C}~(b) and \ref{3D}~(b) tell us in particular
that the $L$-function of the Chow motive $\Gr_0 \Mgm(X)^e$
is computed via (the $e$-part of) interior cohomology of $X$.
\end{Rem}

\begin{Ex} \label{3E}
Let $C$ be a smooth projective curve, and $P \in C$ a $k$-rational point.
Put $X := C - P$. Localization \cite[Prop.~4.1.5]{V} shows that
there is an exact triangle
\[
\Mgm(P) \longto \Mgm(C) \longto \Mcgm(X) \longto \Mgm(P)[1] \; .
\]
The morphism $\Mgm(P) \to \Mgm(C)$ is split; hence $\Mcgm(X)$ is a direct
factor of $\Mgm(C)$. It is therefore a Chow motive. By duality
\cite[Thm.~4.3.7~3]{V}, the same is then true for $\Mgm(X)$.
In particular, both $\Mgm(X)$ and $\Mcgm(X)$ are pure of weight zero.
But the morphism $u: \Mgm(X) \to \Mcgm(X)$ is not an isomorphism
(look at degree $0$ or $-2$, or check that the conclusions of
Theorems~\ref{3C}~(c) and \ref{3D}~(c) do not hold).
Therefore, Assumption~\ref{3A} is not fulfilled for $e = 1$.
Of course, this can be seen directly: the boundary motive
$\dMgm(X)$ has a weight filtration 
\[ 
\BZ(1)[1] \longto \dMgm(X) 
\longto \BZ(0) \longto \BZ(1)[2] 
\]
(which is necessarily split since there are no non-trivial morphisms from
$\BZ(0)$ to $\BZ(1)[2]$). Orthogonality~\ref{1A}~(3) then shows that there
are no non-trivial morphisms from an object of weights $\le -2$ to
$\dMgm(X)$, and no non-trivial morphisms from $\dMgm(X)$ to an object of
weights $\ge 1$. The object $\dMgm(X)$ being non-trivial, we conclude that 
it does not admit a weight filtration avoiding weights $-1$ and $0$. 
\end{Ex}

More generally, we have the following, which again
illustrates just how restrictive Assumption~\ref{3A} is.

\begin{Prop} \label{3Ea} 
Fix data (1) and (2) as before.
Assume that $X$ admits a smooth compactification $\Xp$ such that
the complement $Y = \Xp - X$ is smooth. Then the following
statements are equivalent.
\begin{enumerate}
\item[(i)] Assumption~\ref{3A} is valid, i.e., the object $\dMgm(X)^e$
is without weights $-1$ and $0$.
\item[(ii)] The object $\dMgm(X)^e$ is trivial (hence the conclusion
of Corollary~\ref{3Aa} holds).
\end{enumerate}
\end{Prop}

\begin{Proof}
Statement (ii) clearly implies (i). 
In order to show that it is implied by (i), let us show that
the hypothesis on $\Xp$ and $Y$ forces the boundary motive $\dMgm(X)$
to lie in the intersection
\[
\DeffgM_{w \ge -1} \cap \DeffgM_{w \le 0} \; .
\]
By orthogonality~\ref{1A}~(3), the same is then true for its direct factor
$\dMgm(X)^e$. Thus, the only way for $\dMgm(X)^e$ to avoid 
weights $-1$ and $0$ is to be trivial (again by orthogonality).

In order to show our claim, apply \cite[Prop.~2.4]{W1} to see that  
$\dMgm(X)$ is 
isomorphic to the shift by $[-1]$ of a choice
of cone of the canonical morphism
\[
\Mgm (Y) \oplus \Mgm(X)
\longto \Mgm (\Xp) \; .
\]
In particular, there is a morphism 
$c_+: \dMgm(X) \to \Mgm (Y)$,
and an exact triangle
\[
(W) \quad\quad\quad\quad
C_- \longto \dMgm(X) \stackrel{c_+}{\longto} 
\Mgm (Y) \longto C_- [1] \; ,
\]
where $C_-$ equals the shift by $[-1]$ of a cone of
\[
\Mgm(X) \longto \Mgm ( \Xp) \; .
\]
By assumption, $Y$ is smooth and proper, hence $\Mgm (Y)$,
as a Chow motive,
is pure of weight $0$.
Duality for smooth schemes \cite[Thm.~4.3.7~3]{V}
shows that $C_-$ is pure of weight $-1$.
Hence $(W)$ is a weight filtration of $\dMgm(X)$
by objects of weights $-1$ and $0$.   
Proposition~\ref{1C} then shows that $\dMgm(X)$ belongs indeed to
\[
\DeffgM_{w \ge -1} \cap \DeffgM_{w \le 0} \; .
\]
\end{Proof}

Corollary~\ref{3B} allows to say more about the \'etale realizations.

\begin{Thm} \label{3F}
Keep the situation considered in Theorem~\ref{Main}, 
and fix a prime $\ell > 0$.
Assume that $k$ is the quotient field of a Dedekind domain $A$, and that $k$
is of characteristic zero. Fix 
a non-zero prime ideal $\Fp$ of $A$, and let 
$p$ denote its residue characteristic. \\[0.1cm]
(a)~Assume that $p \ne \ell$. Then a sufficient condition for the $\ell$-adic
realization $H^\ast (\Gr_0 \Mgm(X)^e)$ to be unramified 
at $\Fp$ is the existence of some smooth compactification
of $X$ having good reduction at $\Fp$. A sufficient condition 
for $H^\ast (\Gr_0 \Mgm(X)^e)$ to be semi-stable
at $\Fp$ is the existence of some smooth compactification 
of $X$ having simple semi-stable reduction at $\Fp$. \\[0.1cm]
(b)~Assume that $p = \ell$, and that the residue field $A /\Fp$
is perfect. Then a sufficient condition for the $p$-adic
realization $H^\ast (\Gr_0 \Mgm(X)^e)$ to be crystalline
at $\Fp$ is the existence of some smooth compactification 
of $X$ having good reduction at $\Fp$. A sufficient condition 
for $H^\ast (\Gr_0 \Mgm(X)^e)$ to be semi-stable
at $\Fp$ is the existence of some smooth compactification 
of $X$ having simple semi-stable reduction at $\Fp$.
\end{Thm}

\begin{Proof}
By Corollary~\ref{3B}, $\Gr_0 \Mgm(X)^e$ is a direct factor 
of the motive $\Mgm(\Xp)$ of any smooth compactification $\Xp$ of $X$.
Hence $H^\ast (\Gr_0 \Mgm(X)^e)$ is a direct factor of the cohomology
of any such $\Xp$.

Part~(a) uses the spectral sequence \cite[Scholie~2.5]{D2}
relating cohomo\-logy with coefficients in
the vanishing cycle sheaves $\psi^q$ to cohomo\-logy of the generic fibre
$\Xp$.
By \cite[Kor.~2.25]{RZ} (\cite[Thm.~3.3]{D2} when $p = 0$), 
our assumption on the reduction of $\Xp$ implies that 
the inertia group acts trivially on 
the $\psi^q$. It therefore acts unipotently on the 
cohomology of $\Xp$.
Part~(b) follows from the $C_{st}$-conjecture, proved by
Tsuji \cite[Thm.~0.2]{T} (see also \cite{N} for a proof via $K$-theory).
\end{Proof}

By \cite[Thm.~4.3.7]{V}, the category $\DgM_F$ is a rigid tensor 
triangulated category. Furthermore, if $X$ is smooth of pure dimension $n$,
then the objects $\Mgm(X)$ and $\Mcgm(X)(-n)[-2n]$ are canonically dual
to each other. By \cite[Thm.~6.1]{W1}, the boundary motive $\dMgm(X)$
is canonically dual to $\dMgm(X)(-n)[-(2n-1)]$. Furthermore, these dualities 
fit together to give an identification of the dual of the exact
triangle
\[
(\ast) \quad\quad
\dMgm(X) \longto \Mgm(X) \longto \Mcgm(X) \longto \dMgm(X)[1]
\]
and the exact triangle $(\ast)(-n)[-2n]$. 
The construction of the duality isomorphism \loccit \
shows that the dual of the action of the algebra
$\bar{c}_{1,2}(X,X)$ on $(*)$
equals the natural (anti-)action given by the composition of the canonical
action on $(\ast)(-n)[-2n]$, preceded by the anti-involution ${}^t$.
Consider the idempotent $\, {}^t e$ (i.e., the transposition of $e$).  

\begin{Prop} \label{3Ga}
(a)~Assumption~\ref{3A} is equivalent to any of the following statements.
\begin{enumerate}
\item[(i)] Both $\dMgm(X)^e$ and $\dMgm(X)^{{}^t e}$ are without weight $-1$.
\item[(ii)] Both $\dMgm(X)^e$ and $\dMgm(X)^{{}^t e}$ are without weight $0$.
\end{enumerate}
In particular, Assumption~\ref{3A} is satisfied for $e$
if and only if it
is satisfied for the transposition $\, {}^t e$. \\[0.1cm]
(b)~Assume that $e$ is symmetric, i.e., that $e = {}^t e$. Then 
Assumption~\ref{3A} is equivalent to any of the following statements.
\begin{enumerate}
\item[(i)] $\dMgm(X)^e$ is without weight $-1$.
\item[(ii)] $\dMgm(X)^e$ is without weight $0$.
\end{enumerate}
\end{Prop}

\begin{Proof}
Indeed, $\dMgm(X)^e$ is dual to $\dMgm(X)^{{}^t e}(-n)[-(2n-1)]$.
Now observe 
that $\BZ(-n)[-(2n-1)]$ is pure of weight $1$.
Then use unicity of weight filtrations avoiding weight $-1$ resp.
weight $0$ (Corollary~\ref{2Ca}).
\end{Proof}

\begin{Ex} \label{3Gb}
Assume that an abstract group $G$ acts by automorphisms on $X$. \\[0.1cm]
(a)~The action of $G$ translates into a morphism of algebras 
$\BZ[G] \to c_{1,2}(X,X)$. This morphism transforms the natural 
anti-involution ${}^*$ of $\BZ[G]$ induced by $g \mapsto g^{-1}$ into the
anti-involution ${}^t$ of $c_{1,2}(X,X)$. \\[0.1cm]
(b)~Let $F$ be a flat $\BZ$-algebra, and $e$ an idempotent in $F[G]$.
The morphism of (a) then allows to consider the image of $e$ 
(equally denoted by $e$) in
$c_{1,2}(X,X)$, then in $\bar{c}_{1,2}(X,X)$, and to ask whether 
Assumption~\ref{3A} is valid for $e$. \\[0.1cm]
(c)~As a special case of (b), consider the case when $G$ is finite,
its order $r$ is invertible in $F$, and $e$ is 
the idempotent in $F[G]$ associated to a character $\varepsilon \ $ on $G$
with values in the multiplicative group $F^*$:
\[
e = \frac{1}{r} 
\sum_{g \in G} \varepsilon(g)^{-1} \cdot g \; .
\]
Observe that the idempotent $e \in \bar{c}_{1,2}(X,X)$ is symmetric if 
$\varepsilon^{-1} = \varepsilon$. \\[0.1cm]
(d)~Let us consider the situation from Section~\ref{4}, and show
that Assumption~\ref{3A} is satisfied for
$X = X_n^r \in Sm / \BQ$ and 
\[
e = \frac{1}{(2n^2)^r \cdot r!} 
\sum_{\gamma \in \Gamma_r} \varepsilon(\gamma)^{-1} \cdot \gamma \; .  
\]
As in Section~\ref{4}, denote
by $\bbX_n^{r,\infty}$ the complement of $X_n^r$ in $\bbX_n^r$.
By \cite[Prop.~2.4]{W1}, the object 
$\dMgm(X_n^r)^e$ is canonically
isomorphic to the shift by $[-1]$ of a canonical choice
of cone of the canonical morphism
\[
\Mgm \bigl( \bbX_n^{r,\infty} \bigr)^e \oplus \Mgm(X_n^r)^e
\longto \Mgm \bigl( \bbX_n^r \bigr)^e \; .
\]
In particular, there is a canonical morphism 
$c_+: \dMgm(X_n^r)^e \to \Mgm \bigl( \bbX_n^{r,\infty} \bigr)^e$,
and an exact triangle
\[
C_- \longto \dMgm(X_n^r)^e \stackrel{c_+}{\longto} 
\Mgm \bigl( \bbX_n^{r,\infty} \bigr)^e 
\longto C_- [1] \; ,
\]
where $C_-$ equals the shift by $[-1]$ of a cone of
\[
j_n^r : \Mgm(X_n^r)^e \longto \Mgm \bigl( \bbX_n^r \bigr)^e \; .
\]
By Corollary~\ref{4C}~(b) and (c),
\[
C_- \cong \Mgm \bigl( S_n^\infty \bigr) (r+1)[r+1] 
\]
is pure of weight $-(r+1)$.
It follows from this and from Corollary~\ref{1E}~(a)
that the exact triangle
\[
C_- \longto \dMgm(X_n^r)^e \stackrel{c_+}{\longto} 
\Mgm \bigl( \bbX_n^{r,\infty} \bigr)^e =
\Mcgm \bigl( \bbX_n^{r,\infty} \bigr)^e
\longto C_- [1] 
\]
is a weight filtration of $\dMgm(X_n^r)^e$ avoiding weights 
$-r,\ldots,-1$, and hence in particular, avoiding weight
$-1$ (since $r \ge 1$). Our claim then follows from
Proposition~\ref{3Ga}~(b) (observe that $e$ is symmetric). 

(Alternatively, use Claim~(C) of the proof of Theorem~\ref{4Ca},
to see directly that the last term of the weight filtration of
$\dMgm(X_n^r)^e$,
\[
\Mgm \bigl( \bbX_n^{r,\infty} \bigr)^e \cong
\Mgm \bigl( S_n^\infty \bigr)[r] 
\]
is pure of weight $r \ge 1$.)
\\[0.1cm]
(e)~As a by-product of the above identification of the
weight filtration of $\dMgm(X_n^r)^e$,
we see that the cohomological Betti realization $H^i ( \dMgm(X_n^r)^e )$
\cite{L}, tensored with $\BZ[1 / (2n \cdot r!)]$,
is without torsion for all integers $i$. 
We thus recover a result
of Hida's (see \cite[Prop.~3]{Gh}): the odd primes $p$ dividing
the torsion of the ($e$-part of the) \emph{boundary cohomology} of $X_n^r$
\cite[Sect.~3.2]{Gh}
satisfy $p \le r$ or $p \tei n$.
\end{Ex}

\begin{Rem} \label{3H}
Let us agree to forget the results from Section~\ref{4}, 
and see what the theory developed in the present section 
implies in the situation studied in Example~\ref{3Gb}~(d). 
We only use the validity of Assumption~\ref{3A}. \\[0.1cm]
(a)~It follows formally from Theorem~\ref{Main}~(a)--(c) that 
$\Gr_0 \Mgm(X_n^r)^e$ and $\Gr_0 \Mcgm(X_n^r)^e$ are defined,
and canonically isomorphic, and that there are exact triangles
\[
C_- \longto \Mgm(X_n^r)^e 
\stackrel{\pi_0}{\longto} \Gr_0 \Mgm(X_n^r)^e 
\stackrel{\delta_-}{\longto} C_-[1]
\]
and
\[
\Mgm \bigl( \bbX_n^{r,\infty} \bigr)^e 
\stackrel{\delta_+}{\longto} \Gr_0 \Mcgm(X_n^r)^e 
\stackrel{i_0}{\longto} \Mcgm(X_n^r)^e 
\longto \Mgm \bigl( \bbX_n^{r,\infty} \bigr)^e[1] \; .
\]
(b)~It follows formally from Corollary~\ref{3B} that
$\Gr_0 \Mgm(X_n^r)^e$ is a direct factor of ${}^r_n \CW = \Mgm(\bbX_n^r)^e$,
with a canonical complement. Call this complement $N$.
It follows formally from Theorems~\ref{3C}~(b) and \ref{3D}~(b) 
that the realizations of
the motive 
$\Gr_0 \Mgm \bigl( X_n^r \bigr)^e$ equal interior cohomology,
i.e., the image of the morphism
\[
H^n_c \bigl( X_n^r \bigr)^e \longto 
H^n \bigl( X_n^r \bigr)^e \; .
\]
By \cite[Sect.~1.2.0, Thm.~1.2.1, Sect.~1.3.4]{Scho}, 
the same is true for ${}^r_n \CW$.
Therefore, the complement $N$ has trivial realizations. 
Thus, its underlying Grothendieck motive is trivial.
This means that the construction of the Grothendieck motives 
for modular forms
$M(f)$ can be done replacing the Chow motive ${}^r_n \CW$
by $\Gr_0 \Mgm \bigl( X_n^r \bigr)^e$. \\[0.1cm]
(c)~As recalled in \cite[Sect.~11.5.2]{A}, Voevodsky's 
\emph{nilpotence conjecture} implies that homological
equivalence equals smash-nilpotent equivalence. Therefore
\cite[Cor.~11.5.1.2]{A}, it implies that the forgetful functor
from Chow motives to Grothendieck motives is conservative.
Thus, the nilpotence conjecture gives a hypothetical 
abstract reason for the complement $N$ of
$\Gr_0 \Mgm \bigl( X_n^r \bigr)^e$ in ${}^r_n \CW$ to be zero, 
which would mean that $\Gr_0 \Mgm \bigl( X_n^r \bigr)^e = {}^r_n \CW$. \\[0.1cm]
(d)~Recall that Scholl's motives $M(f)$ are constructed out of ${}^r_n \CW$
using cycles, which are only known to be idempotent after passage to
the Grothendieck motive underlying ${}^r_n \CW$: indeed,
the Eichler--Shimura isomorphism allows for
a control of the action of the Hecke algebra on the relevant cohomology 
group. 
In particular, there are only finitely many eigenvalues, a fact which is
needed for the construction of the projector to the eigenspace corresponding
to $f$. 
The nilpotence conjecture
implies \cite[Cor.~11.5.1.2]{A} that idempotents can be lifted from 
Grothendieck to Chow motives. Hence its validity would mean that Scholl's cycles
can be modified by terms homologically equivalent to zero, to give idempotent
endomorphisms of ${}^r_n \CW$. This would produce Chow motives, whose
underlying Grothendieck motives are the $M(f)$. \\[0.1cm] 
(e)~By \cite[Cor.~7.8]{K}, the conclusion from (d) holds under an assumption,
which is \emph{a priori} 
weaker than the nilpotence conjecture: the \emph{finite dimensionality}
\cite[Def.~3.7]{K} of the Chow motive ${}^r_n \CW$. By \cite[Thm.~4.2, Cor.~5.11]{K},
finite dimensionality is known for motives of curves and motives of Abelian varieties.
Unfortunately, 
the methods of \loccit \ do not seem to admit an immediate generalization to
families of Abelian varieties over curves, or degenerations of such families
(like $\bbX_n^r$).
\end{Rem}

\begin{Rem} \label{3I}
(a)~Of course, none of the implications listed in parts (a)--(c) 
of the preceding
remark are new: they are all consequences of
Theorem~\ref{4Ca}, Corollary~\ref{4C} 
and Corollary~\ref{4Da}, whose proof involves
the geometry
of the boundary of the smooth compactification $\bbX_n^r$
of $X_n^r$. Observe that some of these results
were even used in our proof \ref{3Gb}~(d)
of the validity of Assumption~\ref{3A}.
In other words, we applied the strategy from Remark~\ref{3Ab},
and proved Assumption~\ref{3A} via parts (a) and (c)
of Theorem~\ref{Main}. \\[0.1cm]
(b)~For Shimura varieties of higher dimension, 
Hecke-equivariant smooth compactifications (like $\bbX_n^r$ in the case of 
powers of the universal elliptic curve over a modular curve)
are not known (and maybe not reasonable to expect) to exist. 
In this more general setting, a different strategy of
proof of Assumption~\ref{3A} ``avoiding geometry as far as possible''
would therefore be of interest. 
Such a strategy, disjoint from the one from
Remark~\ref{3Ab}, will be developed in \cite{W3}.
\end{Rem}

%%% Local Variables:
%%% mode: latex
%%% TeX-master: "head"
%%% End:

\bigskip

%%%%%%%%%%%%%%%%%%%%%%%%%%%%%%%%%%%%%%%%%%%%%%%%%%%%%%%%%%%%%%%%%%%%%%%
%
%  Bibliography
%
%%%%%%%%%%%%%%%%%%%%%%%%%%%%%%%%%%%%%%%%%%%%%%%%%%%%%%%%%%%%%%%%%%%%%%%

\end{document}